\documentclass[sn-mathphys-num]{sn-jnl}
\usepackage{graphicx}
\usepackage{todonotes}
\usepackage{amsmath}
\usepackage{amsfonts}
\usepackage{amssymb}
\usepackage{amsthm}
\usepackage{siunitx}
\usepackage{hyperref}
\usepackage{caption}
\usepackage{float}
\usepackage[normalem]{ulem}
\usepackage{cancel}

\DeclareSIUnit{\Osm}{\text{Osm}}

\begin{document}
\title{Evaporation-driven tear film thinning and breakup in two space dimensions}
\author[1]{\fnm{Qinying} \sur{Chen}}

\author[1]{\fnm{Tobin A.} \sur{Driscoll}}

\author[1]{\fnm{R. J.} \sur{Braun}}

\affil[1]{\orgdiv{Department of Mathematical Sciences}, \orgname{University of Delaware}, \orgaddress{\city{Newark}, \postcode{19716}, \state{DE}, \country{USA}}}

\abstract{
Evaporation profiles have a strong effect on tear film thinning and breakup (TBU), a key factor in dry eye disease (DED). In experiments, TBU is typically seen to occur in patterns that locally can be circular (spot), linear (streak), or intermediate . We investigate a two-dimensional (2D) model of localized TBU using a Fourier spectral collocation method to observe  how the evaporation distribution affects the resulting dynamics of tear film thickness and osmolarity, among other variables.   We find that the dynamics are not simply an addition of individual 1D solutions of independent TBU events, and we show how the TBU quantities of interest vary continuously from spots to streaks for the shape of the evaporation distribution. We also find a significant speedup by using a proper orthogonal decomposition to reduce the dimension of the numerical system. The speedup will be especially useful for future applications of the model to inverse problems, allowing the clinical observation at scale of quantities that are thought to be important to DED but not directly measurable in vivo within TBU locales.
}
  
\keywords{Tear film, dry eye disease, lubrication theory, fluorescent imaging, dimension reduction{, proper orthogonal decomposition}}
\maketitle

\section{Introduction}

Every time one blinks, the upper eyelid drops and then rises to paint a thin liquid film, the tear film (TF), on the eye surface \cite{Doane80}.   The healthy TF provides lubrication to the ocular surface and eyelid, antimicrobial defense, a smooth ocular surface for refraction, and a supply of oxygen and nutrition to the avascular corneal epithelium \cite{lempDefinitionClassificationDry2007,willcoxTFOSDEWSII2017}.  When the tear film fails, often called tear breakup (TBU), there are adverse stimuli to the ocular surface \cite{king2018mechanisms}; chronic TBU and problems with the TF are thought to be an etiological factor in the development of
dry eye disease (DED).  DED affects millions of people, with the number depending on the diagnostic criteria used \cite{stapletonDEWSIIepi2017}.  DED diminishes both quality of vision and ocular surface comfort and health \cite{nelsonTFOSDEWSII2017}. The major classes of dry eye include evaporative dry eye (EDE) \cite{lempDefinitionClassificationDry2007}.  EDE is thought to be caused by excessive water loss via evaporation \cite{linDryEyeDisease2014,OCEANreport2013}. This paper seeks to add to the basic science understanding of evaporative water loss from the TF.

A key mechanism of DED are generally accepted to be tear hyperosmolarity resulting from TF instability~\cite{lempDefinitionClassificationDry2007}. Tear osmolarity is the total concentration of osmotically active solutes, primarily salt ions~\cite{stahlOsmolalityTearFilm2012}.   Hyperosmolarity is linked to DED~\cite{tomlinsonTearFilmOsmolarity2006}, but it is not currently possible to measure osmolarity directly in areas of TBU.  Osmolarity is measured clinically in the inferior meniscus or outer canthus \cite{lempTearOsmolarityDiagnosis2011}, but these measurements do not reflect the elevated values that may occur over the cornea \cite{li2Dosmofluor}. 

The tear film is a thin liquid film with multiple layers that establishes itself rapidly after a blink \cite{braunTearFilm2018}. 
The tear film (TF) has three layers: a thin and oily lipid layer that is 20 to 100 nm thick \cite{braunDynamicsFunctionTear2015}, an aqueous layer consisting mainly of water \cite{hollyFormationRuptureTear1973} that is a few microns thick, and a half-micron thick mucin layer called the glycocalyx that is on the ocular surface \cite{king-smithThicknessTearFilm2004}. The lipid layer slows evaporation of water from the TF \cite{mishimaOilyLayerTear1961}, and a healthy glycocalyx facilitates the fluid movements of the ocular surface \cite{gipsonDistributionMucinsOcular2004}. The lacrimal gland supplies the majority of the water from the aqueous layer near the temporal canthus \cite{darttNeuralRegulationLacrimal2009}. Osmosis supplies water from the ocular epithelia \cite{braunDynamicsTear2012}.

TBU happens when a dry spot appears on the eye \cite{nornMICROPUNCTATEFLUORESCEINVITAL1970} and is often evaporation-driven \cite{lempDefinitionClassificationDry2007,willcoxTFOSDEWSII2017}. Tear breakup time (TBUT) is the length of time after a blink and before the appearance of the first dry spot \cite{nornMICROPUNCTATEFLUORESCEINVITAL1970}. A clinical determination of TBUT is based on the observer’s judgment \cite{nornMICROPUNCTATEFLUORESCEINVITAL1970} and may incorporate an average of clinicians’ estimates \cite{choReliabilityTearBreakup1992}. TBUT is a method for determining the stability of the TF and checking for dry eye. A short TBUT is a sign of a poor tear film, and the longer break-up takes, the more stable the tear film \cite{dibajniaTearFilmBreakup2012}.

A number of papers have considered a one-dimensional domain that extended across the open eye with fixed ends which typically called post-blink drainage models \cite{SharmaTiwari98,wong1996deposition,MilPolRad02a,BraunFitt03}.  These models typically imposed large menisci at the domain ends that drained fluid from the relatively flat interior of the TF, and that drove the dynamics.  Braun and Fitt \cite{BraunFitt03} included evaporation and found that it could drive thinning to TBU. Heat transfer and evaporation were considered in Li and Braun \cite{liModelHumanTear2012}, and later treated together with thermal imaging in Dursch et al.\ \cite{DurschFLandThermal2017}. 

TBU models are typically along a short line segment and are designed to be local models of TF flow; these models eliminate any influence from the menisci.
The effects of evaporation and the Marangoni effect on TBU have been studied in a number of papers.  Evaporation drove all the dynamics in the models in \cite{PengEtal2014,braunDynamicsFunctionTear2015,braunTearFilm2018}, though the treatments were different. In \cite{PengEtal2014}, a stationary {lipid layer (LL)} with variable thickness caused increased localized thinning; also included were osmolarity transport in the {aqueous layer (AL)} and osmosis across the AL/cornea interface.  An important result from that work is that diffusion of osmolarity out of the region of evaporation prevents osmosis from stopping thinning as it would in spatially uniform models \cite{braunDynamicsTear2012,braunDynamicsFunctionTear2015}.  Simpler evaporation models were used in \cite{braunDynamicsFunctionTear2015,braunTearFilm2018}, but they included transport of fluorescein so that they could explain what was seen in TF visualization experiments.  Zhong et al.\ \cite{zhongDynamicsFluorescentImaging2019} developed a PDE model with one spatial variable that incorporated both mechanisms.  Including fluorescein dye transport and fluorescence enabled fitting of models to \textsl{in vivo} fluorescence data within TBU regions to estimate parameters that are not possible to directly measure there at the time of writing \cite{lukeParameterEstimation2020,lukeParameterEstimationMixedMechanism2021}.  Recently, ODE models with no space dependence have been successfully fit to fluorescence data in small TBU spots and streaks \cite{lukeFittingSimplifiedModels2021}.  Those models have been coupled to a neural-network based data extraction system to greatly expand the amount of TBU instances that may be studied \cite{driscollFittingODEModels2023}.
Other models have provided fundamental analyses on evaporation \cite{ji2016finitetimerupture,ji2020ssdynnon-conserved} and instability \cite{shi2021instabilitysoliddome,shi2022instabilitybubble}.

Other effects have been included in local models of TBU \cite{ZhangMatar03,ZhangMatar04,PengEtal2014,braunTearFilm2018,zhongMathematicalModelling2018,zhongDynamicsFluorescentImaging2019,deyContinuousMucinProfiles2019,deyContinuousMucinCorrection2020,choudhuryMembraneMucin2021}.
Some models incorporated the effect of soluble mucins in the AL necessitating departure from more standard lubrication models \cite{deyContinuousMucinProfiles2019,deyContinuousMucinCorrection2020}.  Localized non-wettability of the ocular surface was modeled in \cite{choudhuryMembraneMucin2021}.  Combining this or related theories may be applied to in vitro systems of epilthelia \cite{madlmodeltearfilm20,madlmimeticplatform2022} to shed light on the role of dewetting in TF instability.

Imaging of the tear film is an important tool for analyzing its dynamics. Common imaging techniques include fluorescence (FL) imaging \cite{king-smithTearFilmInterferometry2014}, spectral interferometry \cite{king-smithApplicationNovelInterferometric2010} and optical coherence tomography \cite{wangPrecornealPrePostlens2003}. Injection of dyes such as fluorescein have been used to stain epithelial cells \cite{nornMICROPUNCTATEFLUORESCEINVITAL1970}, estimate tear drainage rates or turnover times \cite{webberContinuousFluorophotometricMethod1986}, visualize general TF dynamics \cite{benedettoVivoObservationTear1984,begleyQuantitativeAnalysisTear2013}, estimate TF breakup times \cite{nornMICROPUNCTATEFLUORESCEINVITAL1970}, and capture the TBU regions. Simultaneous imaging can help interpret TF dynamics \cite{himebaughScaleSpatialDistribution2012}.



\begin{figure}
    \centering
    \includegraphics[width=\textwidth]{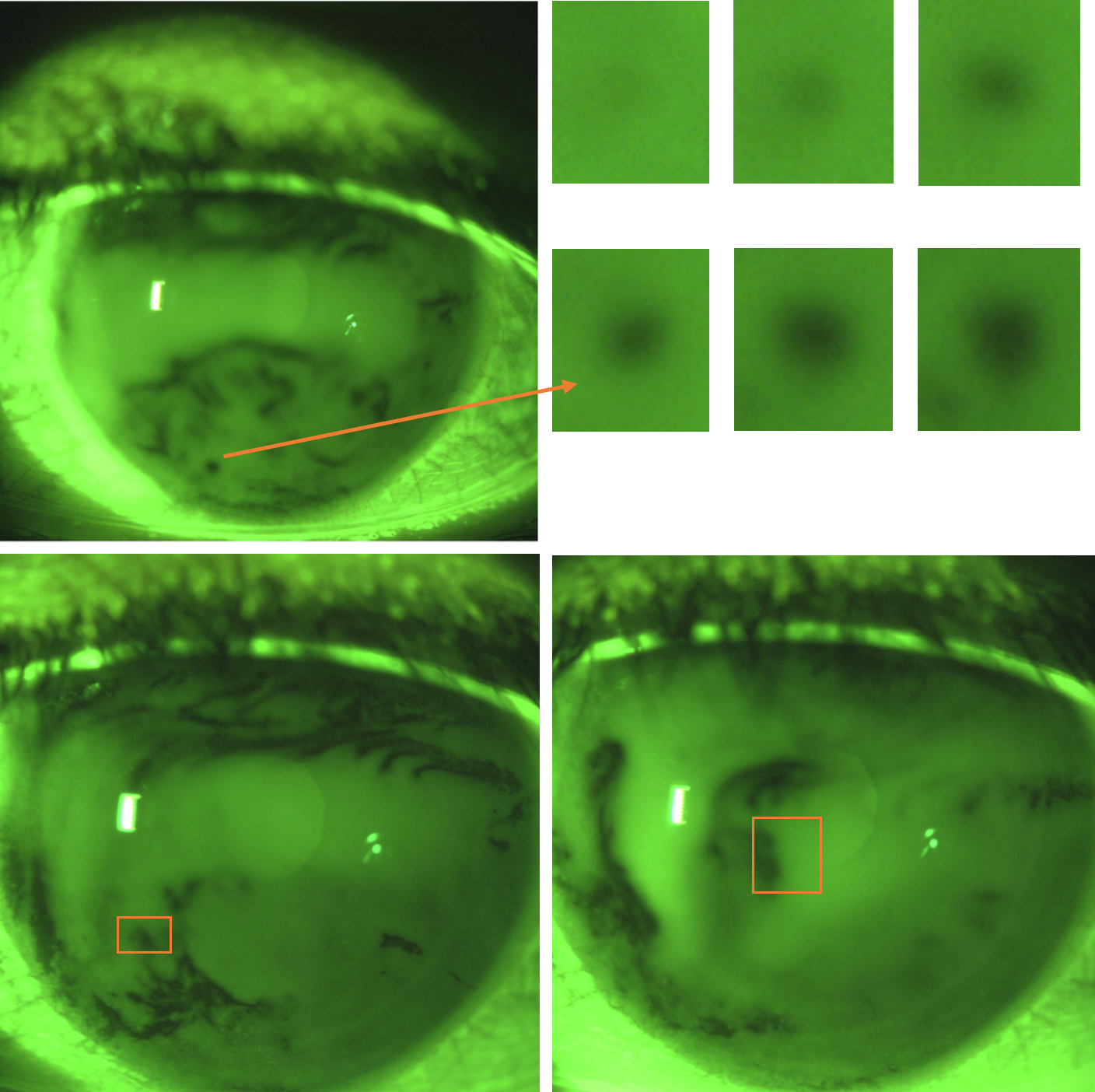}
    \caption{{At top left, bottom left and bottom right are fluorescence images from three different video capture of STARE trials ~\cite{awisi-gyauChangesCornealDetection2019a}. The top right shows a time series of the evolution of the dark spot highlighted in the top left image. The highlights in the bottom left and bottom right images show respectively an indeterminate shape and a formation of connected spots.}(Image courtesy of D.Awisi-Gyau.)}
    \label{fig:flexample}
\end{figure}

\autoref{fig:flexample} {shows the FL images from video captures of three different sustained tear exposure (STARE) trials \cite{awisi-gyauChangesCornealDetection2019a}. As the eye was held open, we observe the emergence and development of small, dark regions. There are mixed patterns for the dark areas: spot-like patterns, which have roughly radial symmetry, streak-like patterns, which are locally one-dimensional, and intermediate shapes. The top right shows a time series of the evolution of a dark spot highlighted in the top left image. The red box in the bottom left image shows a dark spot of indeterminate shape, and the box in the bottom right image shows two connected spots. Our interest is to study these emerging small dark areas from the FL imaging.}  

Luke et al. \cite{lukeParameterEstimation2020} developed techniques for fitting experimental data from FL imaging to  models given in Braun et al. \cite{braunTearFilm2018} for evaporation-driven tear film thinning. They found realistic optimal values for peak and baseline evaporation rates and dry spot sizes, and their thinning rate estimates fell within experimental ranges \cite{nicholsThinningRatePrecorneal2005a}.  { While the fitting of PDE models in this way \cite{lukeParameterEstimation2020} was successful, it was slow work, and only twenty instances of TF thinning were fit in that paper.  By simplifying the model and automating the TBU instance selection, Driscoll et al. \cite{driscollFittingODEModels2023} fit hundreds of thinning instances for normal subjects.  This approach greatly extended the number of thinning instances that could be analyzed, and the trends remained consistent with PDE models.  However, the details of spatial dynamics of thinning and TBU were lost. One clear lesson from these works is that faster computation and fitting of the imaging data could produce much more information about tear film dynamics in human subjects.  Another important point is that fitting the models to data could determine important in vivo quantities, such as thinning rate, that often require much more complex models to estimate local in vivo values (e.g., \cite{PengEtal2014,stapf2017duplex,DurschFLandThermal2017}).  Faster simulations to analyze more TBU instances is major motivator for studying the numerical methods applied in this work.}

Over the decades, researchers have developed various numerical methods to solve thin film problems. Zhornitskaya \& Bertozzi \cite{zhornitskayaPositivityPreservingNumericalSchemes1999} developed  positivity-preserving finite difference schemes in one and two space dimensions.  Oron \& Bankoff \cite{oronDynamicsCondensingLiquid2001} used a ﬁnite difference scheme based on a spatially uniform staggered grid, and a modified Euler method in time. Schwartz et al. \cite{schwartzDewettingPatternsDrying2001} employed a time-marching finite difference method with an alternating-direction-implicit (ADI)
technique. Lee et al. \cite{leeEfficientAdaptiveMultigrid2007} used a full approximation storage (FAS) adaptive multigrid algorithm. Heryudono et al. \cite{heryudonoSingleequationModels2007} used a modified Chebyshev spectral collocation for the spatial discretization. Gr\"un \& Rumpf \cite{grunNonnegativityPreservingConvergent2000} presented finite element methods for fourth order degenerate parabolic equations.  
These methods are all successful for finding a single solution efficiently, but may be slow for repeated solution inside an algorithm.

To our knowledge, Maki et al. \cite{makiTearFilm2010a} were the first to model the human tear film dynamics on a 2D domain. They formulated a relaxation model on an eye-shaped domain with flux boundary conditions,  which was solved numerically with { a} composite overset grid method of 2nd-order finite differences in space together with an adaptive 2nd order BDF time stepping method. The simulation took about two hours on a high-end desktop computer. Li et al. improved the model { first} by incorporating evaporation and ocular surface wettability \cite{liTearFilmDynamics2014}  { and again by} incorporating osmolarity transport and osmosis on the 2D eye-shaped domain \cite{liComputedTearFilm2016}. The spatial derivatives were  discretized using 2nd order finite difference composite overset grids in space { as in \cite{makiTearFilm2010a}, together with} a hybrid time-stepping method, with Runge-Kutta-Chebyshev (RKC)  for solving the osmolarity transport  { and} an adaptive BDF method for updating thickness and pressure. The simulation took about 1 to 2 hours on a high-end desktop computer.  These 2D simulations are computationally expensive if one wants to do repeated { simulations to investigate the effects of the many important parameters or solve inverse problems to fit collected data.}  {proper orthogonal decomposition (POD)} method to reduce the computation time as much as possible.

In this paper, we extend the simulation of evaporative models of TBU dynamics to two dimensions{ (2D) for what we believe to be the first time. As a first step, we simplify the situation to a local model of TF thinning and TBU, as has been done in lower dimensions~\cite{PengEtal2014,braunDynamicsFunctionTear2015,braunTearFilm2018,lukeParameterEstimation2020}. Ignoring boundary (e.g., eyelid) effects allows us to employ a Fourier spectral method that can finish accurate simulations in tens of minutes. We further show that the proper orthogonal decomposition (POD) method~\cite{grasslePODReducedOrder2020,weissTutorialProperOrthogonal2019} can reduce the dimensionality of the discretized equations at the cost of a tolerable amount of error, thus greatly accelerating the simulations. Such accelerations will be valuable not only in manually exploring the dynamics within the model but also for solving the inverse problem in future fits to experimental data in the manner of Driscoll et al.~\cite{driscollFittingODEModels2023}. As a step toward understanding TBU dynamics for human ocular surface health, we demonstrate the use of our solution method on exploratory scenarios, such as the transition from (radially symmetric) spot to (mirror-symmetric) streak, the breakdown of radial symmetry, and the merging of nearby dark spots, that provide initial insights into 2D TBU dynamics.}

The paper is organized as follows. We present the mathematical model in \autoref{sec:model} in nondimensional form. We present the numerical method, including the use of the POD, in \autoref{sec:numerics}. In \autoref{sec:results} we present simulation results for a variety of scenarios, including replication of the previous results for spots and streaks, and give evidence for the effectiveness of the POD. Finally, in \autoref{sec:conclusion} we discuss the results and describe the next steps for future work.  

\section{Model}
\label{sec:model}

Our model is a two-dimensional version of the models derived by Braun et al. \cite{braunTearFilm2018} and used 
by Luke et al. \cite{lukeParameterEstimation2020} for local TBU dynamics.  Li et al \cite{li2Dosmofluor} derived the 2D model equations in detail, though they were not applied to local TBU dynamics. 

The TF is modeled as a Newtonian fluid over a flat corneal surface at the plane $z=0$. 
{ The flat cornea approximation has been justified elsewhere \cite{berger1974,BraunUsha12}, and a 1D model using a curved cornea showed little difference in TF dynamics save for small differences near the lids \cite{Allouche17}.}  Because the film is thin, there is a separation of scales, and lubrication theory may be applied \cite{ODB97,CrasMat09review}. 

The key variables in the system are shown in \autoref{tab:variables}. These quantities are nondimensionalized according to 
\begin{align*}
x'&=\ell x , \quad y'=\ell y , \quad z'=dz, \quad t'=\frac{d}{v_\text{max}}t, \quad h'=dh,\quad u'=\frac{v_\text{max}}{\epsilon}u. \\
v'&=\frac{v_\text{max}}{\epsilon}v, \quad w'=v_\text{max}w, \quad J'=\rho v_\text{max}J, \quad c'=c_0 c,\quad f'=f_{\text{cr}}f,
\end{align*}
where primes denote dimensional quantities and the relevant physical parameters are given in \autoref{tab:physical}. { Other needed nondimensional parameters are given in \autoref{tab:parameters}.}  

{ The derivation of the system follows a typical pattern of deriving thin film models \cite{ODB97,CrasMat09review}.  Scaling produces many factors of $\epsilon=d/\ell \ll 1$ preceding terms that are neglected at leading order. One finds easily that the pressure $p$ is independent of depth and is set by capillarity via the normal stress condition so that it is given by to the curvature of the film surface:
\begin{equation}
p =-\partial_x^2 h - \partial_y^2 h. \label{eq:pdep1}
\end{equation}
Integrating mass conservation for the incompressible case (the divergence-free condition) over $0 < z < h(x,y,t)$ and using the boundary conditions results in the mass conservation equation for the film thickness $h$:
\begin{equation}
\partial_t h  + \partial_x (h\overline{u})+\partial_y (h\overline{v}) = -J + P_c(c-1) \label{eq:pdeh1}
\end{equation}
Here $\overline{u}$ and $\overline{v}$ are depth-averaged velocity components along the film in the $x$ and $y$ directions, respectively, arising from solving an approximation of the momentum conservation equation in each direction:
\begin{align}
\overline{u} =-\frac{h^2}{12}\partial_x p, \label{eq:pdeu}\\
\overline{v} =-\frac{h^2}{12}\partial_y p. \label{eq:pdev}
\end{align}
The terms involving $\overline{u}$ and $\overline{v}$ determine the change in thickness due to motion of fluid inside the film.
$J$ is the thinning rate of the film due to evaporation (discussed below), which appears in the kinematic condition at the air/film interface \cite{ODB97,braunDynamicsTear2012}.  
$P_c(c-1)$ is the osmotic supply of water specified by the boundary condition on the transverse velocity component at the film/cornea interface, which approximates the film/cornea interface as a semipermeable membrane that allows only water to cross \cite{braunDynamicsTear2012}.  This condition represents the supply of fluid from the cornea as the osmolarity in the film rises as a response to evaporation of water from the film.  }

{ Following the approach of Jensen and Grotberg \cite{JenGrot93} for solutes inside the film (and used in e.g., \cite{braunDynamicsTear2012,braunTearFilm2018}) one obtains transport equations for the osmolarity $c$ and the FL concentration $f$.
In the resulting equation for the osmolarity, we can interpret certain collections of terms as being responsible for particular mechanisms: 
\begin{equation}
    \partial_t c +
    \underbrace{\overline{u}\partial_x c ,+ \overline{v}\partial_y c}_{\text{advection}} 
    =
    \underbrace{\frac{\partial_x (h\partial_x c) +\partial_y (h\partial_y c)}{h \text{Pe}_c}}_{\text{diffusion}} + 
    \underbrace{\frac{Jc}{h}}_{\text{evap.}} - 
    \underbrace{\frac{P_c(c-1)c}{h}}_{\text{osmosis}}. \label{eq:termtypes}
\end{equation}
The transport equation for $f$ is very similar.  The evaporation term increases solute concentrations, while osmosis lowers solute concentrations (when osmolarity rises)\cite{li2Dosmofluor}.}

The resulting nondimensional system { to be solved for the dependent variables $h$, $p$, $c$, and $f$} is 
\begin{align}
\partial_t h  + \partial_x (h\overline{u})+\partial_y (h\overline{v})&= -J + P_c(c-1), \label{eq:pdeh}\\
p &=-\partial_x^2 h - \partial_y^2 h, \label{eq:pdep}\\
h(\partial_t c +\overline{u}\partial_x c + \overline{v}\partial_y c) &= \text{Pe}_c^{-1}(\partial_x (h\partial_x c)+\partial_y (h\partial_y c))+Jc-P_c(c-1)c, \label{eq:pdec}\\
h(\partial_t f +\overline{u}\partial_x f + \overline{v}\partial_y f) &= \text{Pe}_f^{-1}(\partial_x (h\partial_x f)+\partial_y (h\partial_y f))+Jf-P_c(c-1)f, \label{eq:pdef}
\end{align}
{ with} parameters given in \autoref{tab:parameters}. The {quantities} $\overline{u}$ and $\overline{v}$ are { determined by \eqref{eq:pdeu} and \eqref{eq:pdev}, but they are not dependent variables per se. Note that $f$ appears only in~\eqref{eq:pdef}, so that equation can be solved separately once the other dependent variables are known.}

\begin{table}[tbp]
    \centering
    \begin{tabular}{cl}
        Variable & Meaning \\ \hline
        $x$, $y$ & transverse spatial dimensions \\
        $z$ & depth dimension \\
        $t$ & time \\
        $h(x,y,t)$ & TF thickness \\ 
        $u(x,y,t)$, $v(x,y,t)$ & transverse fluid velocities \\ 
        $p(x,y,t)$ & pressure \\
        $J(x,y)$ & evaporation rate \\ 
        $c(x,y,t)$ & osmolarity \\ 
        $f(x,y,t)$ & fluorescein concentration
    \end{tabular}
    \caption{Variables in the two-dimensional model.}
    \label{tab:variables}
\end{table}

\begin{table}[tbp]
    \centering
    \begin{tabular}{c l c l} 
      Parameter & Description & Value & Reference \\ 
     \hline
     $\mu$ & Viscosity & \qty{1.3e-3}{\Pa\sec} & Tiffany \cite{tiffanyViscosityHumanTears1991} \\
     $\sigma_{0}$ & Surface tension & \qty{0.045}{\N\per\meter} & Nagyov\'a and Tiffany \cite{nagyovaComponentsResponsibleSurface1999} \\
     $\rho$& Density & \qty{e3}{\kg\per \meter\tothe{3}} & Water \\
     $d$ & Initial TF thickness & \qty{4.5}{\micro\meter} & Calculated \\  
     $\ell$ & $(\sigma_0/ \mu / v_{\max})^{1/4} d$ & \qty{0.54}{\milli \meter} & Calculated \\ 
     $v_{\text{max}}$ & Peak thinning rate & \qty{10}{\micro\meter\per\minute} & Nichols et al. \cite{nicholsThinningRatePrecorneal2005a} \\ 
     $V_{w}$ & Molar volume of water & \qty{1.8e-5}{\meter\cubed \per \mol} & Water \\ 
     $D_{f}$ & Diffusivity of fluorescein & \qty{0.39e-9}{\meter\squared \per \sec} & Casalini et al. \cite{casaliniDiffusionAggregationSodium2011} \\ 
     $D_{o}$& Diffusivity of salt & \qty{1.6e-9}{\meter\squared \per \sec} & Riquelme et al. \cite{riquelmeInterferometricMeasurementDiffusion2007} \\
     $c_{0}$& Isotonic osmolarity & \qty{300}{\Osm \per \meter\cubed} & Lemp et al. \cite{lempTearOsmolarityDiagnosis2011}\\
     $P_{0}$& Permeability of cornea & \qty{12.1}{\um\per\sec} & Braun et al. \cite{braunDynamicsFunctionTear2015}\\
     $\epsilon_{f}$& Napierian extinction coefficient & \qty{1.75e7}{\L \per \meter \per \mol} & Mota et al. \cite{motaSpectrophotometricAnalysisSodium1991}\\
     $f_{cr}$& Critical FL concentration  &  $0.2\%$ & Webber and Jones
     \cite{webberContinuousFluorophotometricMethod1986}
    \end{tabular}
    \caption{Physical parameters (dimensional) used in the governing equations.}
    \label{tab:physical}
\end{table}

\begin{table}
    \centering
    \begin{tabular}{c c c} 
     Parameter & Expression & Value  \\ \hline
     $\epsilon$ & $d/\ell$ & $8.3\times 10^{-3}$ \\[1mm]
     $\alpha$ & $({\alpha_{0} \mu})/({\rho \ell \epsilon^3})$ & $4.06\times 10^{-2}$  \\[1mm]
     $P_{c}$ & $({P_{0} V_{w} c_{0}})/({v_\text{max}})$ & $0.392$  \\[1mm]
      $\text{Pe}_{f}$ & $({v_\text{max} \ell})/({\epsilon D_{f}})$ & $27.7$ \\[1mm] 
      $\text{Pe}_{c}$ & $({v_\text{max} \ell})/({\epsilon D_{0}})$ & $6.76$ \\[1mm]
     $\phi$ & $\epsilon_{f} f_{\text{cr}} d$ & $0.417$ 
    \end{tabular}
    \caption{Typical values of nondimensional parameters that appear in the model \eqref{eq:pdeh}--\eqref{eq:pdef} using the parameters in \autoref{tab:physical}.}
    \label{tab:parameters}
\end{table}

We  { apply the equations over a small area} of the cornea that is not close to the eyelids and limbus. Because we are not interested in the effects of these boundaries, we assume periodic conditions on all the dependent variables. We also assume that the simulation begins after the eye opens { (e.g., as in Figure~6 of \cite{King-SmithIOVS13a})} and that all the dependent variables are initially uniform:
\begin{equation}
\label{eq:ic}
    h(x,y,0) = c(x,y,0) = 1,\quad f(x,y,0) = f_0,
\end{equation}
where $f_0$ is the FL concentration normalized to the critical concentration $f_{\text{cr}}$. { For the numerical DAE solver, it is also helpful to supply the initial pressure $p(x,y,0)=0.$}

The evaporation rate $J$ is our primary input to the model, and it drives all the dynamics. Inhomogeneities in an in vivo lipid layer, sitting atop the aqueous layer of the tear film, are presumed to cause local increases in the evaporation rate, leading to local decreases in $h$ and corresponding increases in the solute concentrations.  { An example of this sequence occurring in vivo is shown in Figure~6 of King-Smith et al. \cite{King-SmithIOVS13a}, and simplified models using this idea have appeared, e.g.\ \cite{PengEtal2014,braunDynamicsFunctionTear2015,braunTearFilm2018}.} We represent the spatial variation of $J$ as one or more localized peaks: 
\begin{equation} \label{eq:Jm}
J(x,y) = v_b + \sum_{k=1}^{K} (a_k-v_b)\,G\left( \frac{x-x_k}{x_{w,k}}, \frac{y-y_k}{y_{w,k}} \right),
\end{equation}
where $v_b = v_{\text{min}}/v_{\text{max}}$ is a baseline value, $(x_k,y_k)$ is the center of the $k$th peak, $a_k>v_b$ is the height of the $k$th peak, $x_{w,k}$ and $y_{w,k}$ are characteristic widths of peak $k$, and $G$ is the Gaussian
\begin{equation}
    \label{eq:gaussian}
    G(x,y) = \exp\left[-(x^2 + y^2) / 2 \right].
\end{equation}
{ We note that there was little difference between the resulting images whether one used an evaporation distribution with a very smooth Gaussian or with a sharp-edged evaporation distribution using a hyperbolic tangent function \cite{braunTearFilm2018}.  We choose to use the smooth Gaussian distribution in this work.  Our chosen peak and background values for evaporation are well within experimental observations \cite{nicholsThinningRatePrecorneal2005a}.}


The FL intensity $I$ is obtained via \cite{webberContinuousFluorophotometricMethod1986,braunModelTearFilm2014}
\begin{equation}
\label{eq:FLI}
I=I_0\frac{1-\exp(-\phi fh)}{1+f^2}, 
\end{equation}
where $I_0$ is a normalization coefficient and $\phi$ is the nondimensional Napierian extinction coefficient in \autoref{tab:parameters}.

Given an evaporation function $J(x,y)$, we solve the system \eqref{eq:pdeh}--\eqref{eq:pdef} to obtain $h$, $p$, $c$, and $f$ as functions of space and time. These can then be inserted into~\eqref{eq:FLI} to find fluorescent intensity. Note that while equation \eqref{eq:pdef} describes how the evolution of FL concentration $f$ depends on $h$, $p$, and $c$, those quantities do not in turn depend on $f$. Hence in practice we first solve for the dynamics of $h$, $p$, and $c$, and then use those in a separate stage to solve for $f$.  Once $h$ and $f$ are known, we compute the FL intensity via \eqref{eq:FLI}.

\section{Numerical methods}
\label{sec:numerics}

For the dynamics we use the method of lines with a Fourier spectral collocation method~ \cite{trefethenSpectralMethodsMatlab2000} in space on a uniform periodic grid on the domain $(\pi, \pi]^2$. The number of grid points $m$ and $n$ in each dimension is chosen to be even. The resulting discretization of spatial terms in \eqref{eq:pdeh}--\eqref{eq:pdef} creates a differential--algebraic system (DAE) for  that is solved in Julia using the QNDF solver, an adaptive quasi-constant time step stiff method in the \textbf{DifferentialEquations} package \cite{rackauckas2017differentialequations} { similar to backward differentiation formulas and using Shampine's accuracy-optimal kappa values as defaults \cite{shampine1997matlab}}. In cases where the evaporation function $J$ has a single peak centered at the origin, we exploit symmetry by solving only over the first quadrant $[0,\pi]^2$. 

The full DAE system has size $\approx 4mn$, with a factor of 4 reduction if symmetry is used. A simulation such as those in \autoref{sec:results} takes several minutes on a typical workstation. In the future, we want to use the DAE solution as the forward method in an inverse solver for potentially hundreds of experimental trials. Thus, acceleration of the forward solver is desirable. 

We can greatly reduce the dimension of the discretization via the proper orthogonal decomposition (POD) method~\cite{grasslePODReducedOrder2020,weissTutorialProperOrthogonal2019}.
The full-size discrete system may be written generically as
\begin{subequations}
    \label{eq:pod1}
    \begin{align}
        \partial_t h(t) &= { F}_h(t,h,p,c), \label{eq:pod1h} \\ 
        0 &= { F}_p(t,h,p,c), \label{eq:pod1p} \\
        \partial_t c(t) &= { F}_c(t,h,p,c), \label{eq:pod1c}
    \end{align}
\end{subequations}
where we now use $h$, $p$, and $c$, to denote discretizations of length $mn$. Suppose that we find the solution over a short time interval $t =(0, \tau)$. Choosing times { $t_1=0 < t_2< \dots < t_N=\tau$}, we can form the $mn\times N$ matrix
\begin{equation}
\label{eq:podh}
    H = \begin{bmatrix}
    h(t_1) & h(t_2) & \dots & h(t_N)
    \end{bmatrix}.
\end{equation}
The POD exploits the fact that the matrix $H$ is often well-represented by a low-rank projection via the SVD. Intuitively, the structure of $h(x,y,t)$ tends to stay close to a manifold of dimension far less than $mn$. 

Thus, we let $B_h$ be the first $k$ left singular vectors of $H$, where $k \ll N$. Now $B_{h}^{T}B_{h}$ is the $k\times k$ identity and $B_{h} B_{h}^{T}$ is an orthogonal projector onto the column space of $B_{h}$.  Let $\Tilde{h}=B_{h}^{T} h$. We now have approximate dynamics to replace \eqref{eq:pod1h}:
\begin{align}
    \partial_t h &\approx { F}_h(t,B_h \Tilde{h},p,c). \label{eq:pod2}
\end{align}
The system \eqref{eq:pod2}, however, has large dimension $mn$. It is reasonable to project it down into the $B_h$ basis via
\begin{align*}
    B_h^T \partial_t h  &\approx B_h^T { F}_h(t,B_h \Tilde{h},p,c),\\
    \partial_t \Tilde{h} &\approx B_h^T { F}_h(t,B_h \Tilde{h},p,c).
\end{align*}

We can similarly repeat the projection process for $p$ and $c$, finding bases $B_p$ and $B_c$ and replacing~\eqref{eq:pod1} with the low-dimensional approximation
\begin{subequations}
    \label{eq:pod3}
    \begin{align}
        \partial_t \Tilde{h}(t) &= B_h^T { F}_h(t,B_h \Tilde{h},B_p \Tilde{p},B_c\Tilde{c}), \label{eq:pod3h} \\
        0 &= B_p^T { F}_p(t,B_h \Tilde{h},B_p \Tilde{p},B_c\Tilde{c}), \label{eq:pod3p} \\
        \partial_t \Tilde{c}(t) &= B_c^T { F}_c(t,B_h \Tilde{h},B_p \Tilde{p},B_c\Tilde{c}), \label{eq:pod3c}
    \end{align}
\end{subequations}
with subsequent lifting to the original discretization grid occurring via
\begin{equation}
    \label{eq:podlift}
    h \approx B_h \Tilde{h},\quad p \approx B_p \Tilde{p},\quad c \approx B_c \Tilde{c}.
\end{equation}
Note that the dimensions of the bases $B_h$, $B_p$, and $B_c$ do not need to be identical. The savings for using~\eqref{eq:pod3} in place of~\eqref{eq:pod1} are in the time integration only; the reduced vectors have to be reconstituted onto the original discretization via~\eqref{eq:podlift} and manipulated at that size in order to compute ${ F}_h$, ${ F}_p$, and ${ F}_c$. Even so, the savings for the DAE solver are significant.

The typical use of POD is to solve the full system over a relatively short interval $[0,\tau]$, then project to low dimension to continue the evolution. However, we are able to realize substantial additional gains in cases where the evaporation function~\eqref{eq:Jm} comprises only one peak or $K$ well-separated peaks. In these cases we can find the full solution for each individual peak at the much lower cost of 1D radial solutions~\cite{lukeParameterEstimation2020}. The tear film equations are in \autoref{secA1} for this case. We can solve \eqref{eq:A1}--\eqref{eq:A5} for each peak while centered at the origin, then shift the center and interpolate the solution to the Cartesian grid as necessary. In the case of the $h$ component, for example, we use 1D solutions to obtain snapshot matrices $H_1,\ldots,H_K$ as in \eqref{eq:podh}, then concatenate them horizontally before taking the SVD. This process takes negligible time compared to solving~\eqref{eq:pod3}, but it produces unsatisfactory approximation bases when the peaks overlap sufficiently, as we show in \autoref{sec:results}.  

\section{Results}
\label{sec:results}

{ We halt simulations when the dimensional thickness of the thinning tear film reaches \qty{1}{\micro\m}, since the model is unlikely to be valid as the thickness decreases further. We denote this event as TBU break-up time (TBUT), although we do not claim that it represents a true rupture time in a physically valid model. The results presented in this section were obtained using Julia version $1.8.5$ on a Windows laptop with a 12th Gen Intel(R) Core (TM) i7-12700H 2.30GHz processor and 32GB RAM.}  

\begin{figure}[tbp]
    \centering
    \includegraphics[width=0.95\textwidth]{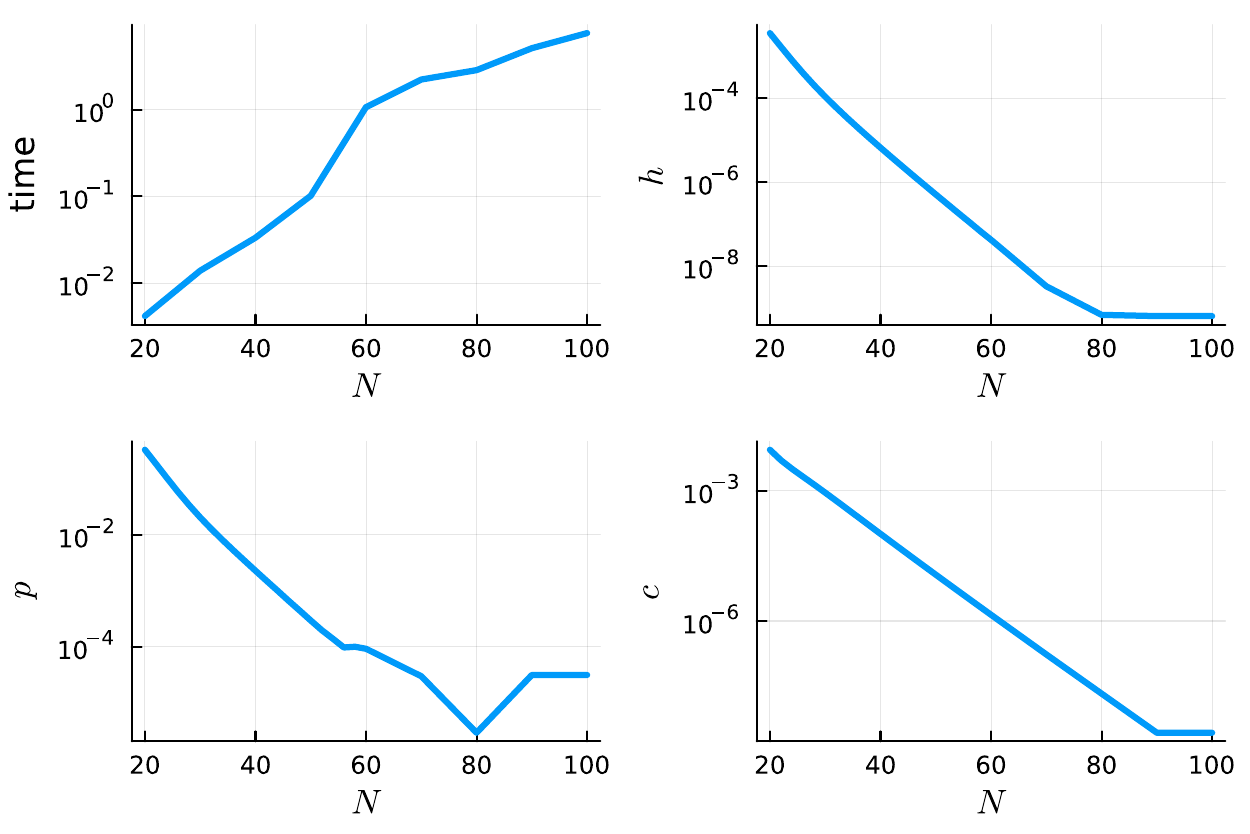}
    \caption{{ Top left: Computation time in hours with respect to the Fourier grid resolution $N$. The other three figures show the relative error at TBUT in the dependent variables $h,$ $p,$ and $c$.} }
    \label{fig:griderr}
\end{figure}

{ \autoref{fig:griderr} shows the computation time and relative error of the solution for increasing grid resolutions. Here we solve \eqref{eq:pdeh}--\eqref{eq:pdef} for a single spot centered at origin, with evaporation widths $x_w=y_w=0.5$, baseline value $v_b=0.1$, and  { peak value} $a_1=1$ in \eqref{eq:Jm}. The error is computed at breakup time $t=2.4$ using a reference solution at grid size $N=100$ and Fourier interpolation onto matching nodes.  We observe that for $N=60$, the solution is essentially converged in space, suggesting that time-stepping errors become dominant. Since the computation at $N=100$ takes more than 7 times longer than at $N=60$, we fix the spatial grid to be $60\times 60$ for all simulations presented in what follows.}

\subsection{Single TBU region}

We start with $K=1$ peak in the evaporation function~\eqref{eq:Jm}. In all the cases of this section, the center of the peak occurs at the origin, i.e., $x_1=y_1=0$. 

\subsubsection{Circular spot}

For the circularly symmetric case with $x_w = y_w$, the problem reduces to the 1D axisymmetric model described in \autoref{secA1} and solved previously \cite{braunTearFilm2018,lukeParameterEstimation2020}. This fact allowed us to check the full 2D results against independently generated results from the simpler model that are interpolated onto the 2D grid. 

\begin{figure}[tbp]
    \centering
    \includegraphics[width=0.95\textwidth]{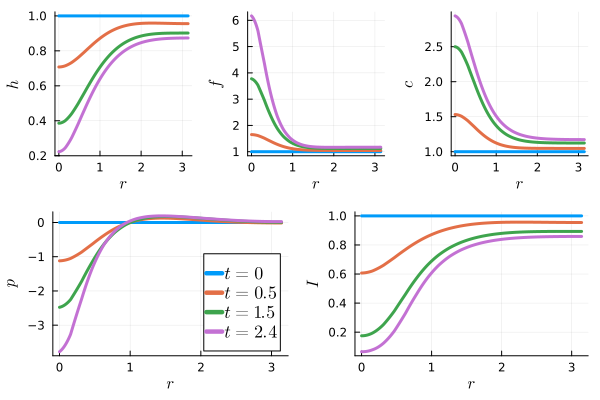}
    \caption{Nondimensional solution for a single circular evaporation spot, shown as functions of distance $r$ from the center of the spot and time $t$. The peak evaporation rate is $a_1=1$ over a baseline rate of $v_b=0.1$, and the spot width is $x_w=y_w=0.5$.}
    \label{fig:single_spot_over_space}
\end{figure}

The results for a circular evaporation spot with width $x_w=y_w=0.5$, peak value $a_1=1$, and baseline value $v_b=0.1$ are shown in \autoref{fig:single_spot_over_space} as functions of the radius. The integration was halted at TBUT $t\approx 2.4$. We observe that the film is thinnest at the center, which is the same location as the peak evaporation rate. The pressure is lowest at the center and is maximum at around $r=1$, and the pressure difference, caused by capillarity, drives fluid in towards the middle. The solute concentrations are maximum at the center and increase as time evolves. Finally, FL intensity $I$ is lowest where the thickness $h$ is smallest and FL concentration $f$ is highest. 

\begin{figure}[tbp]
    \centering
    \includegraphics[width=0.95\textwidth]{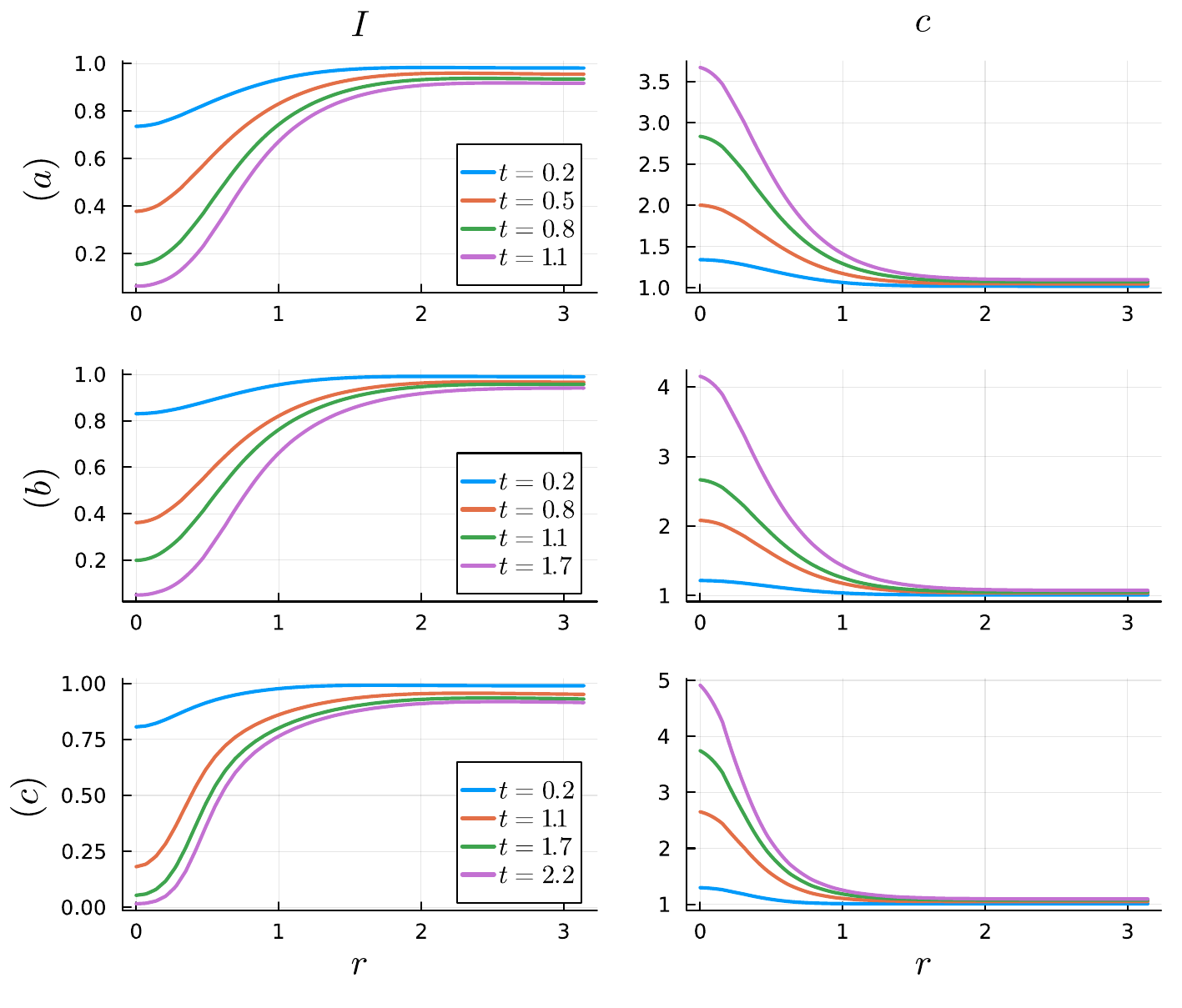}
    \caption{{FL intensity and osmolarity at different time levels for three different one-spot solutions. The TBUTs for  cases (a), (b), and (c) in \eqref{eq:threespots} are $1.1$, $1.7$, and $2.2$, respectively.}}
    \label{fig:three_spot_over_r}
\end{figure}

In order to explore the effects of parameters in the evaporation function $J$, in \autoref{fig:three_spot_over_r} we show the intensity { and osmolarity function at different time levels for three different cases over radius $r$} :
\begin{subequations}
    \label{eq:threespots}
    \begin{align}
        a_1=1.5, &\quad v_b=0.1, \quad x_w=y_w=0.5, \label{eq:threespots_a}\\
        a_1=1.0, &\quad v_b=0.05, \quad x_w=y_w=0.5, \label{eq:threespots_b}\\
        a_1=1.5, &\quad v_b=0.05, \quad x_w=y_w=0.3. \label{eq:threespots_c}
    \end{align}
\end{subequations}
{ The TBUTs for each of these three simulations are $1.1$, $1.7$, and $2.2$, respectively. In case~\eqref{eq:threespots_a} (first row in the figure), we see a more rapid thinning at the center of the spot due to a larger peak evaporation rate. In case~\eqref{eq:threespots_b} (second row in the figure), intensity drops quicker, and osmolarity becomes higher over $r$ at the the same time level compared to case~\eqref{eq:threespots_c} (third row). }


\begin{figure}[tbp]
    \centering
    \includegraphics[width=0.9\textwidth]{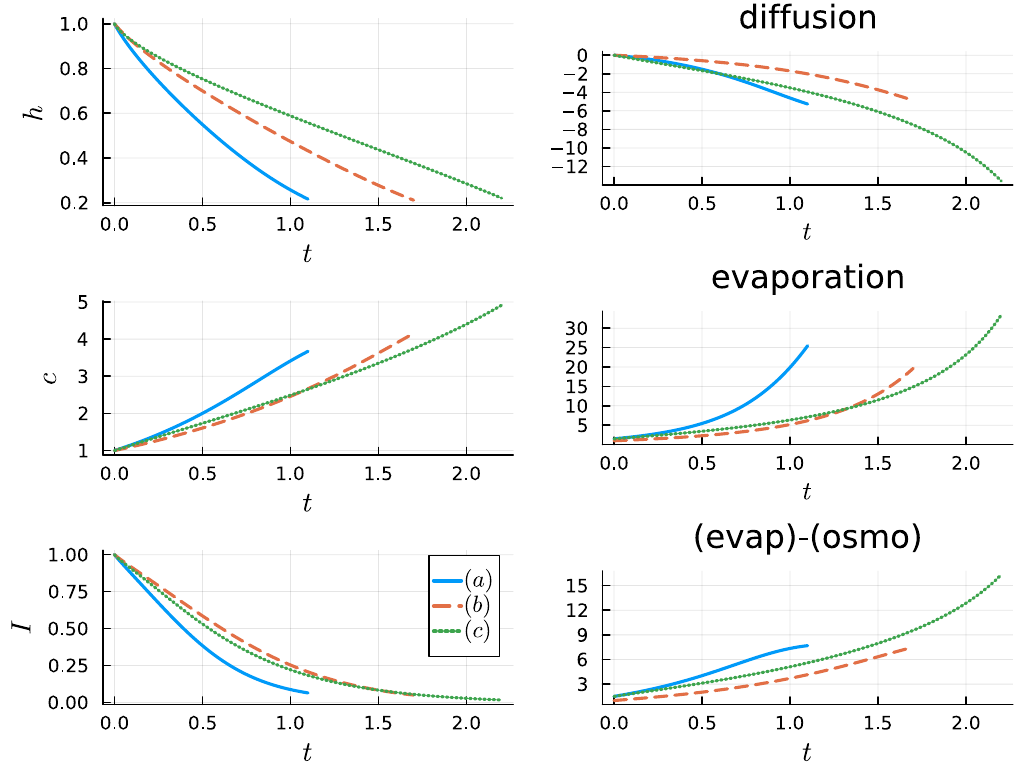}
        \caption{Left: Central values of tear film thickness, osmolarity and FL intensity for three different spots \eqref{eq:threespots}. Right: Dominant terms in the osmolarity equation \eqref{eq:termtypes}.}
    \label{fig:three_spot_diff}
\end{figure}

\autoref{fig:three_spot_diff} shows the central values of dependent variables and the mechanism interpretations given in~\eqref{eq:termtypes} for diffusion, evaporation and osmosis. We observe that the TBUT is the smallest for case {(a); the higher evaporation rate in case (c) is more than counteracted by increased capillary flow.} The osmolarity becomes large due to evaporation that is faster than the diffusion. For case (a), around $t=0.4$, the FL intensity drops to half of its initial values; for cases (b) and (c), the intensity drops to this level at around $t=0.6$.    

\clearpage

\subsubsection{Transition from spot to streak}

 { Previous models assumed either axisymmetric spots or streaks with mirror symmetry about a line} (e.g., \cite{ZhangMatar03,PengEtal2014,braunDynamicsFunctionTear2015,braunTearFilm2018,zhongDynamicsFluorescentImaging2019,lukeParameterEstimation2020}), but none, to our knowledge, { have studied cases that are intermediate to these.} 
The 2D model allows us to simulate intermediate states that are observed frequently in data \cite{liu2009link,wingli2015TBUinfrared,awisi-gyauChangesCornealDetection2019a}. Our 2D periodic boundary conditions do not allow fully 1D streaks, so instead we use different values for $x_{w,1}$ and $y_{w,1}$ in~\eqref{eq:Jm} to make ellipses in the level curves of $J$, and the streak limit corresponds to $x_{w,1}\to 0$. We use peak evaporation rate $a_1=1$ and the baseline rate $v_b=0.1$ throughout this section. We approach the spot-to-streak transition in two different ways. 

 { \textbf{Fixed total evaporation rate.} In this approach,} the product $x_w\, y_w$ is kept fixed at $\tfrac{1}{4}$ in order to keep the total evaporation rate  constant; { i.e., the total amount of water lost, as represented by the double integral of $J$, is constant for different shapes.} 

\begin{figure}[tbp]
    \centering
    \includegraphics[width=0.95\textwidth]{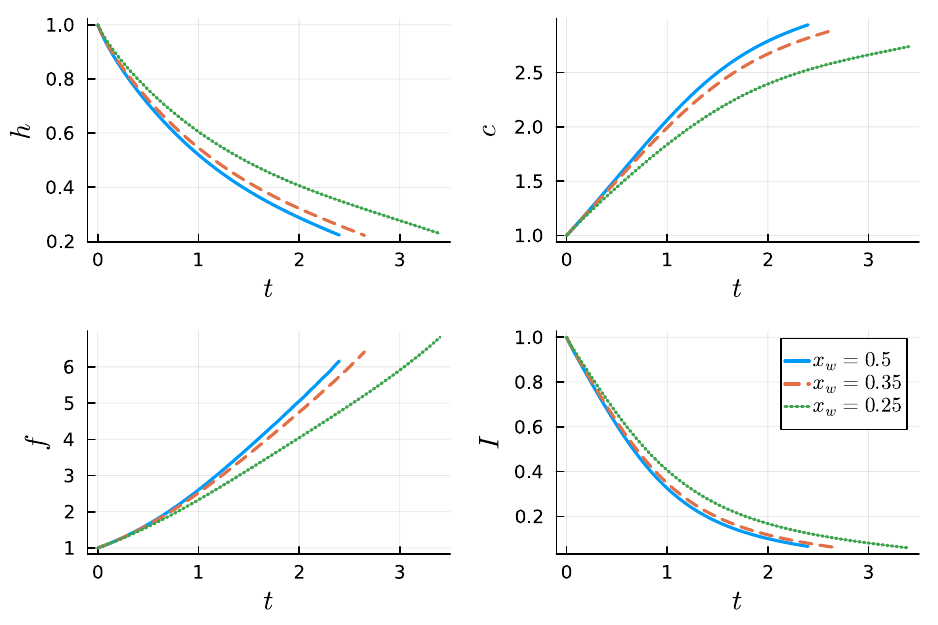}
        \caption{Central values of $h$, $c$, $f$, and $I$ for spot-to-streak transition with fixed product of peak widths $x_w y_w=\tfrac{1}{4}$, baseline rate $v_b=0.1$, and peak evaporation rate $a_1=1$. The evaporation width $x_w$ is $0.5,0.35,0.25$, and the TBUT are $2.4,2.65,3.4$ respectively.}
    \label{fig:fixed_product_hcfi_center}
\end{figure}

\begin{figure}[tbp]
    \centering
    \includegraphics[width=0.95\textwidth]{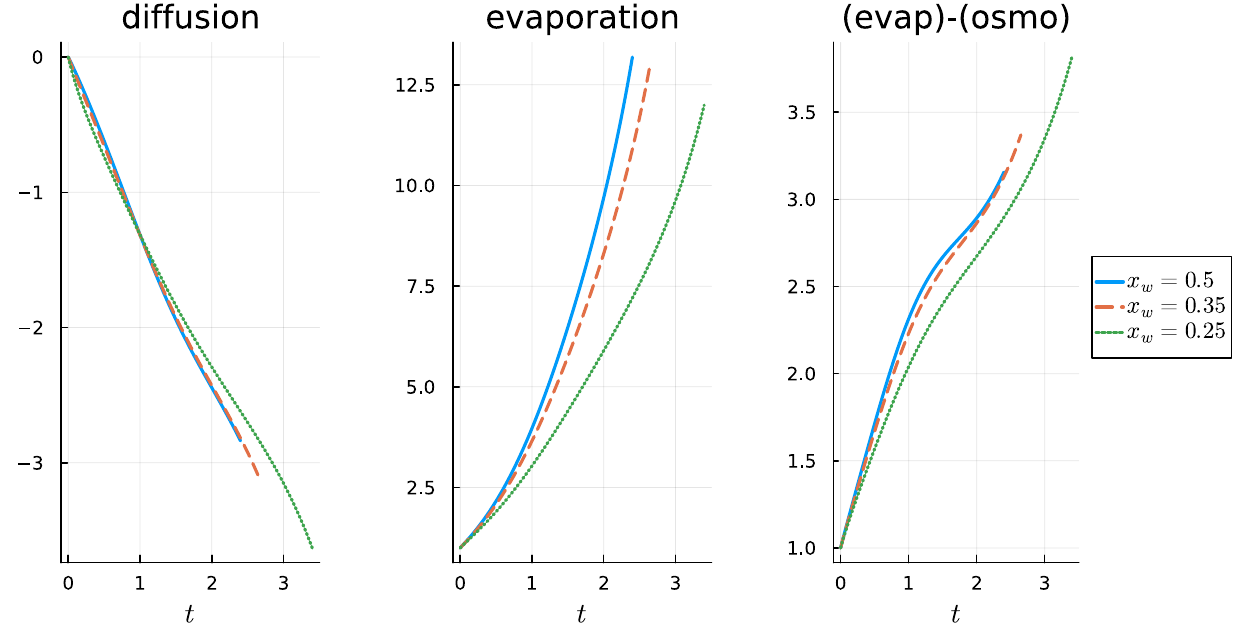}
    \caption{Dominant terms in the osmolarity equation \eqref{eq:termtypes}. The parameters for the evaporation function are the same as in the simulations for \autoref{fig:fixed_product_hcfi_center}.}
    \label{fig:fixed_product_diff}
\end{figure}

\begin{figure}[tbp]
    \centering
    \includegraphics[width=0.95\textwidth]{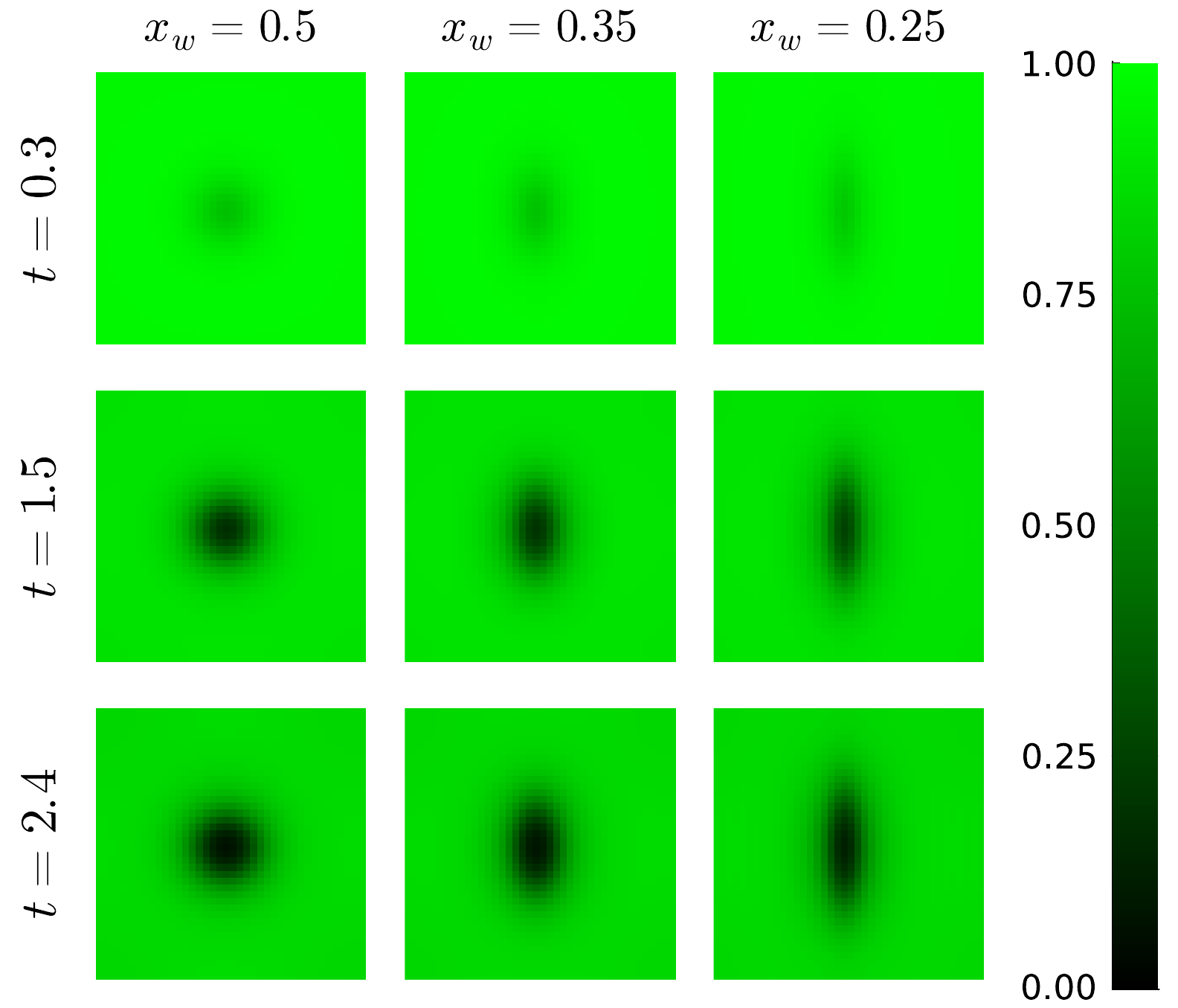}
    \caption{FL intensity for spot-to-streak transition under the same parameters as in \autoref{fig:fixed_product_hcfi_center} with fixed product $x_w y_w=\tfrac{1}{4}$.}
    \label{fig:sptosk_fixed_product_intensity}
\end{figure}

\autoref{fig:fixed_product_hcfi_center} shows the central values of variables as functions of time for $x_{w}=0.5$, $0.35$, and $0.25$ for increasing eccentricity.  The TBUT for these cases are $t=2.4,$ $2.65,$ and $3.4$, respectively, showing that the breakup occurs more slowly for streaks than for spots. The axisymmetric spot shows both a higher FL concentration and osmolarity at the center at a fixed time than in the eccentric cases. \autoref{fig:fixed_product_diff} compares the sizes of the mechanisms labeled in~\eqref{eq:termtypes}, exclusive of advection, which is close to zero here. At a fixed time, the spot TBU has a less diffusion, more evaporation, and more net thinning from the difference between evaporation and osmosis, than the eccentric cases. At TBUT, that is, at the ends of the curves, the ordering of the terms is reversed.  
\autoref{fig:sptosk_fixed_product_intensity} shows snapshots of the FL intensity in these three cases.  We note that the levels of intensity, once they become dark enough (say $I\le 0.25$), are difficult to distinguish visually.

\begin{figure}[tbp]
    \centering
    \includegraphics[width=0.95\textwidth]{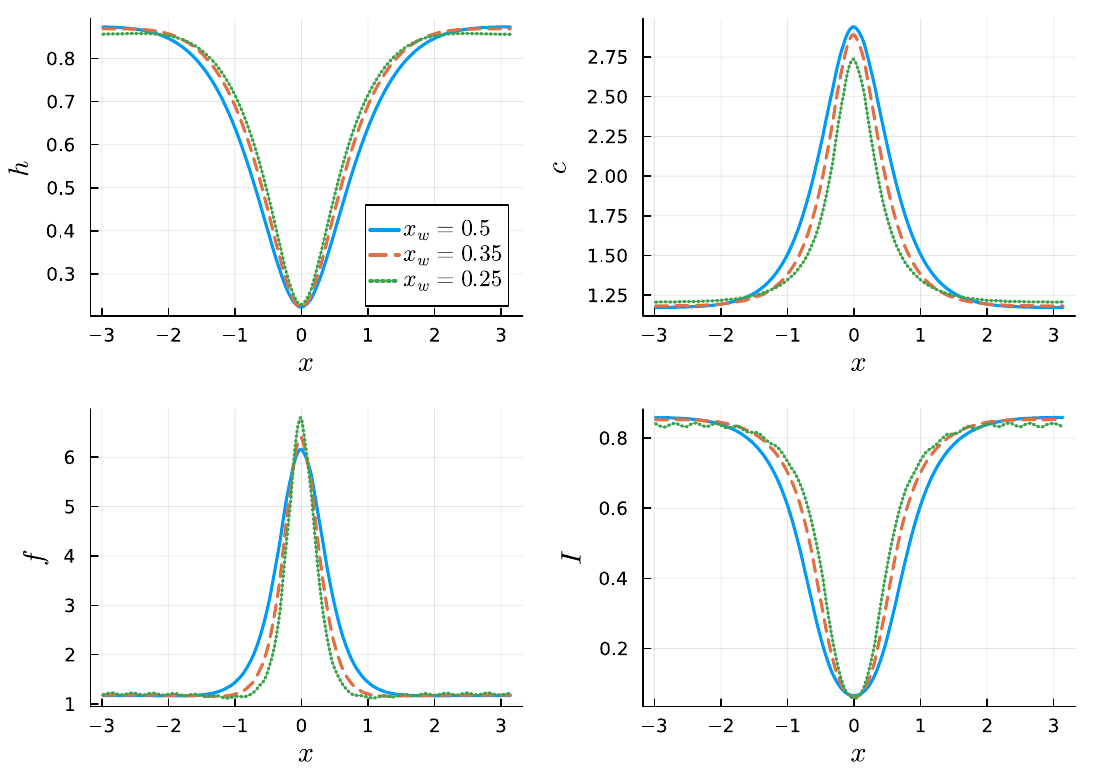}
    \caption{Values of $h$, $c$, $f$, and $I$ along the $x$ axis at TBU time under the same parameters as in \autoref{fig:fixed_product_hcfi_center}.}
    \label{fig:sptosk_fixed_product_x_axis}
\end{figure}

In order to more easily compare these three cases, we turn to one-dimensional plots of the results at the end of the computations.
\autoref{fig:sptosk_fixed_product_x_axis} shows the values of $h$, $c,$ $f,$ and $I$ in these cases along the $x$-axis at their respective TBUTs. The { curves for} $h$  show only weak influence from the evaporative width, as its length scale is dominated by the surface tension.  The effect of diffusion is apparent in the { lowering of the} peak osmolarity value { as the simulation lasts longer for larger TBUT}  The FL concentration is more strongly affected by advection than osmolarity, and this is evident from the fact that higher values occur centrally for larger TBUTs. The intensity profiles depend on both $h$ and $f$, so that results for this variable are intermediate. 

\begin{figure}[tbp]
    \centering
    \includegraphics[width=0.95\textwidth]{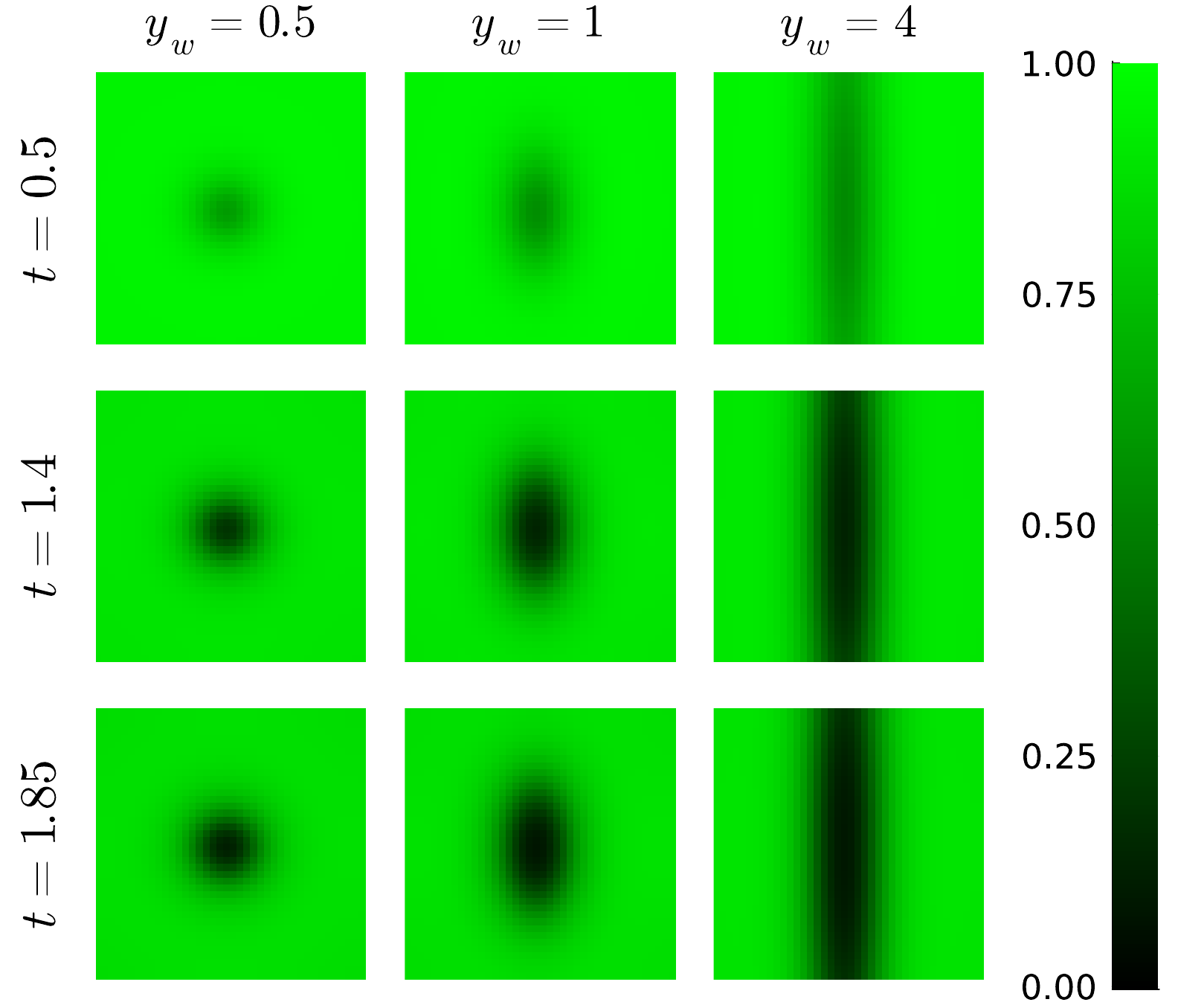}
    \caption{FL intensity for spot and ellipse TBU with fixed $x_w=0.5$ at three time levels. From left to right, the evaporation width $y_w$ is $0.5,1,4$ respectively.}
    \label{fig:sptosk_fixed_xw_intensity}
\end{figure}

\begin{figure}[tbp]
    \centering
    \includegraphics[width=0.95\textwidth]{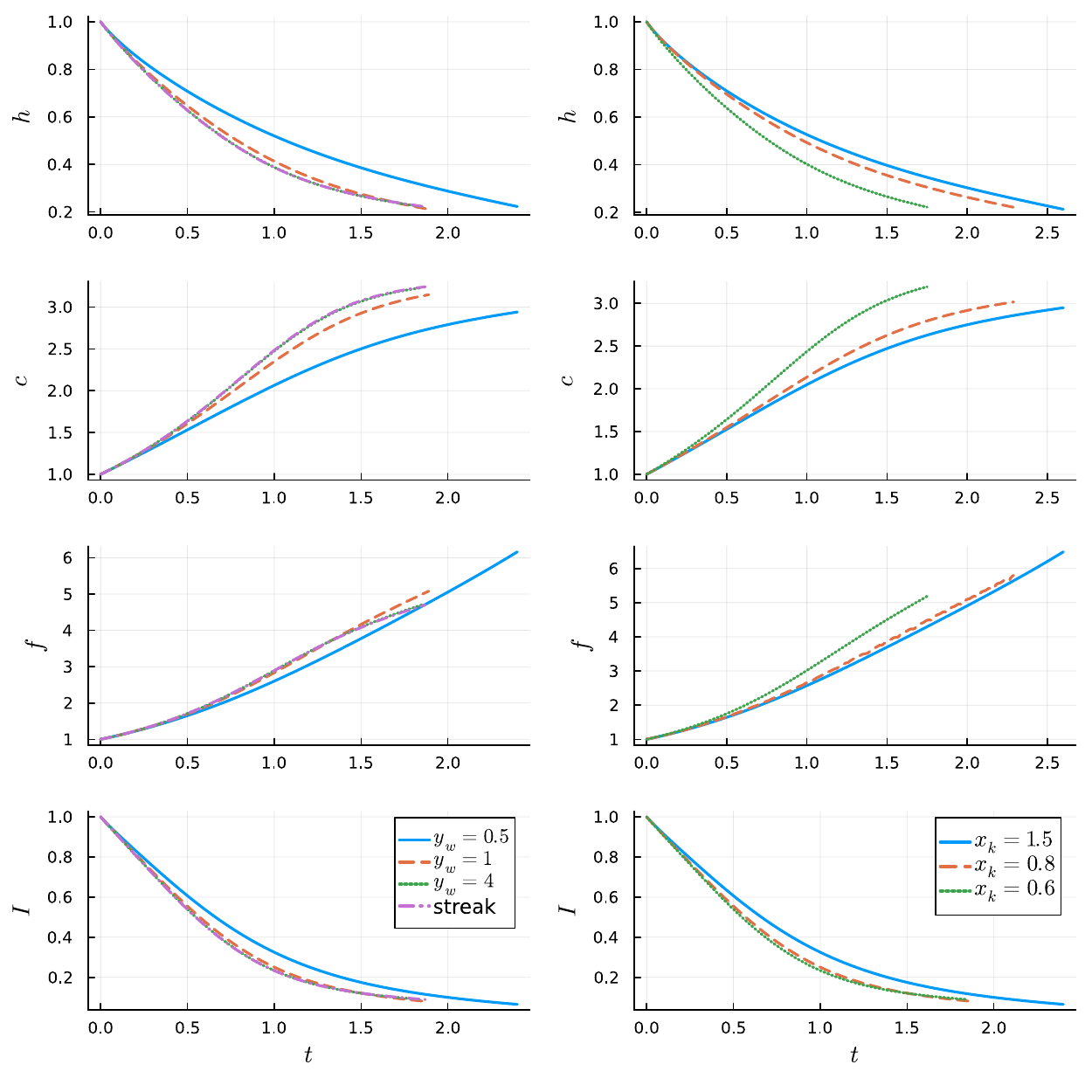}
    \caption{{ Left column: Central values of key quantities for TBU with evaporation functions in elliptical patterns of varying eccentricity. Also shown is a 1D streak solution, which is indistinguishable from the most eccentric case. Right column: Central values for circular spots of differing separation half-distances of $x_k = 1.5$, $0.8$, and $0.6$; see \autoref{sec:multiplespots}.}}
    \label{fig:sptosk_fixed_xw_hcfi_center_and_connected_spots_center}
\end{figure}

\clearpage
{ \textbf{Fixed width in the evaporation function.} In this approach, we leave the Gaussian width $x_w$ fixed at $\tfrac{1}{2}$ while $y_w\to \infty$, along with $v_{max}=1$ and $v_b=0.1$. The second approach may be relevant to tears because the $v_{max}$ and $v_b$ may be roughly constant in a given eye, and so different shapes may allow different amounts of water to evaporate. It also
} allows us to compare directly with the fully 1D streak limit~\cite{braunTearFilm2018}. Note that in this case,  { $y_w \to \infty$} increases the volume of water lost to evaporation for the same time interval compared to fixing { the product} $x_w y_w$.  \autoref{fig:sptosk_fixed_xw_intensity} shows snapshots of the FL intensity for $y_w = 0.5,1,4$, in which the TBUT are $2.4,1.9,1.85$, respectively. We observe that at $y_w=4$, the solution looks indistinguishable from a streak, which has TBUT $1.87$. {The left column in} \autoref{fig:sptosk_fixed_xw_hcfi_center_and_connected_spots_center} shows the central values of key quantities for these cases, along with the 1D streak limit, again showing little difference between $y_w=4$ and a true streak. We also observe that the streak has the highest osmolarity at TBUT, while the spot has the highest FL concentration at TBUT.  {The right column shows results for connected spots, which we will discuss in \autoref{sec:multiplespots}.}

\begin{figure}[tbp]
    \centering
    \includegraphics[width=0.95\textwidth]{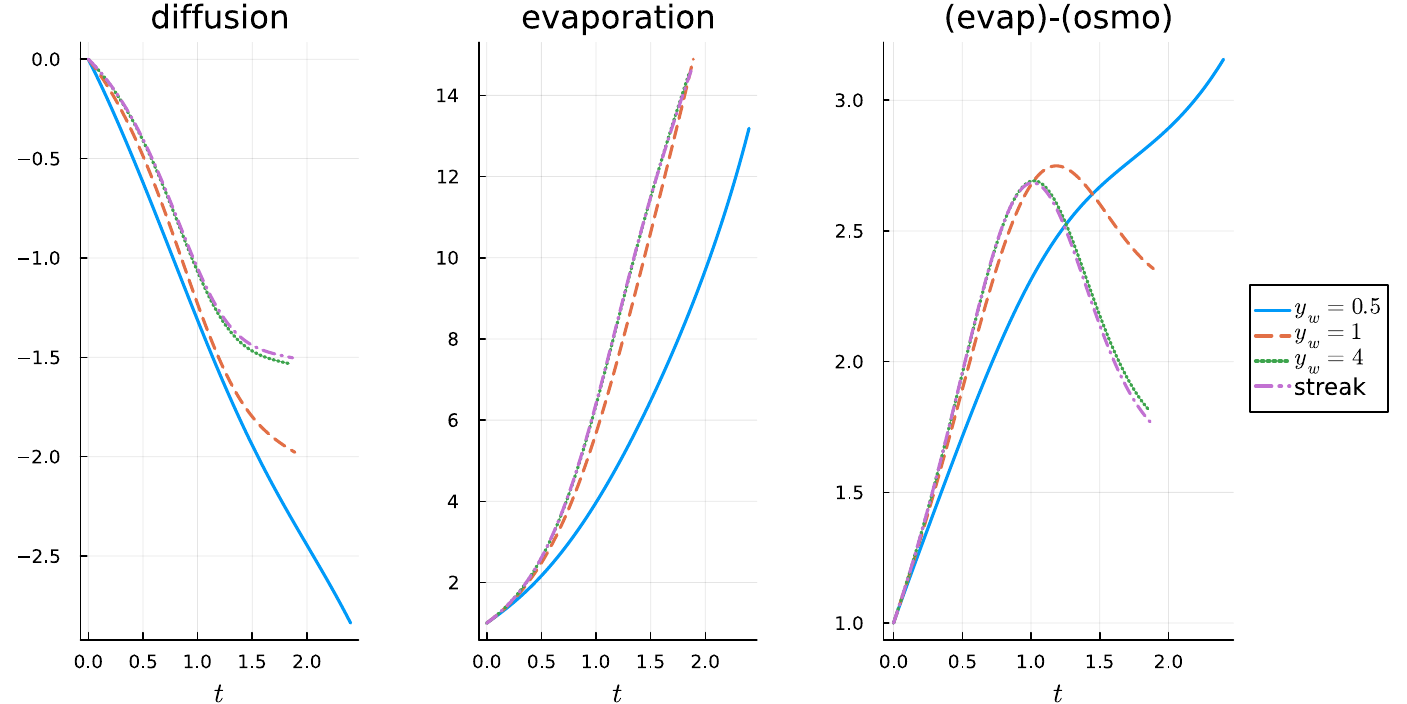}
    \caption{Dominant terms in the osmolarity equation \eqref{eq:termtypes}. The parameters for the evaporation function are the same as in the simulations for {the left column in} \autoref{fig:sptosk_fixed_xw_hcfi_center_and_connected_spots_center}.}
    \label{fig:sptosk_fix_xw_diff}
\end{figure}

\autoref{fig:sptosk_fix_xw_diff} shows the dominant terms in the osmolarity equation \eqref{eq:termtypes} for the second type of spot-to-streak transition. The spot has a larger diffusion and smaller evaporation compared to the other cases, leading to a lower osmolarity. We also observe that for a streak, the difference between evaporation and osmosis reaches a peak and then decreases as we approach TBU. This is because the larger evaporation results in a higher osmolarity, leading to increased osmotic flow.

\subsubsection{Effectiveness of the POD method for a single spot}

\begin{figure}
    \centering
    \includegraphics[width=0.95\textwidth]{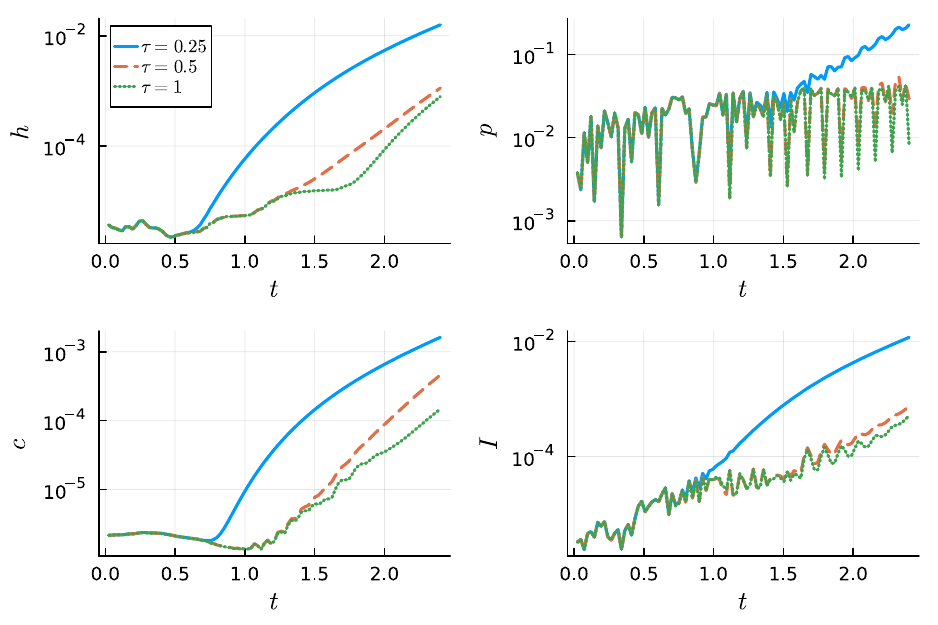}
        \caption{{Relative error of $h,p,c,I$ between POD solution and true solution for a single spot with respect to time. Color indicates different $\tau$ values. Log scale is used for the $y$ axis.} }
    \label{fig:diff_tau}
\end{figure}

In \autoref{fig:diff_tau} we present the relative error between our POD model and true solution for a single spot centered at the origin with evaporation width $x_w = y_w = 0.5$, peak evaporation rate $a_1 = 1$ and baseline rate $v_b=0.1$. The result in the plot is obtained by first solving~\eqref{eq:pod1} over the time domain $t=[0,{\tau}]$ to obtain the POD basis as described in \autoref{sec:numerics}, then solving the reduced-dimension system starting again from $t=0$ to obtain the solution up to TBUT. The resolution for the spatial domain is $60\times 60$ over $[-\pi,\pi] \times [-\pi,\pi] $. { We plot for the three cases $\tau = 0.25,0.5,1$, and the number $N$ of time steps in the captured snapshots in the solution matrix $H$ in \eqref{eq:podh} is $40,50,100$ respectively. Increasing $N$ beyond these values does not improve the solution noticeably. Since the pressure is essentially a second derivative of $h$, it is the component with the largest error. The numbers of basis vectors for $B_h$, $B_p$, $B_c$ for the three cases in the figure are $(15,25,15)$, $(20,30,20)$, and $(40,50,40)$ respectively. We do not show the errors of the $f$ component because they are similar to those for $c$.}

{
In \autoref{fig:diff_tau} we see significant improvement in the accuracy when $\tau$ is changed from $0.25$ to $0.5$. But there is much less improvement when $\tau$ is increased again to 1. In most cases the error in the solution components begins to increase sharply at some time $t > \tau$, although this sudden change does not appear in the derived $I$ quantity. Based on these observations, and because the runtime is more or less proportional to $\tau$, we use $\tau=0.5$ in all the simulation results we present. Thus, when the TBUT is around 2, for instance, we have achieved a 4x speedup.}  

For circularly symmetric evaporation patterns, without going through ~\eqref{eq:pod1}, we take advantage of the 1D radial solver, and we solve on $t=(0,3)$ since it only takes  { about two} more seconds for 1D solver compared to solving on $t=(0,0.5)$.

\clearpage

 On a $50\times 50$ spatial grid, the computation time for a single spot centered at origin using the original Fourier discretization~\eqref{eq:pod1} over $t\in [0,2.4]$ is approximately $2$ minutes. The computation time for using~\eqref{eq:pod1} on $t=(0,0.5)$ is about 30 seconds, and then computing the POD basis and solving the reduced system~\eqref{eq:pod3} takes about another 30 seconds on $t\in [0,2.4]$; thus, about half of the total computation time was saved by using the POD. As the discretization dimension increases, however, the cost of solving~\eqref{eq:pod1} increases dramatically. On a $60\times 60$ spatial grid, for instance, the full solution takes about $40$ minutes, compared to $10$ minutes over $t\in [0,0.5]$ followed by just $35$ seconds for the POD method. In other words, the cost of the POD is essentially entirely due to the full solution step for the initial part of the time evolution, and the low-dimensional projection over the rest of the time domain comes for free.

\clearpage

\subsection{Multiple TBU regions}
\label{sec:multiplespots}

\begin{figure}
    \centering
    \includegraphics[width=0.95\textwidth]{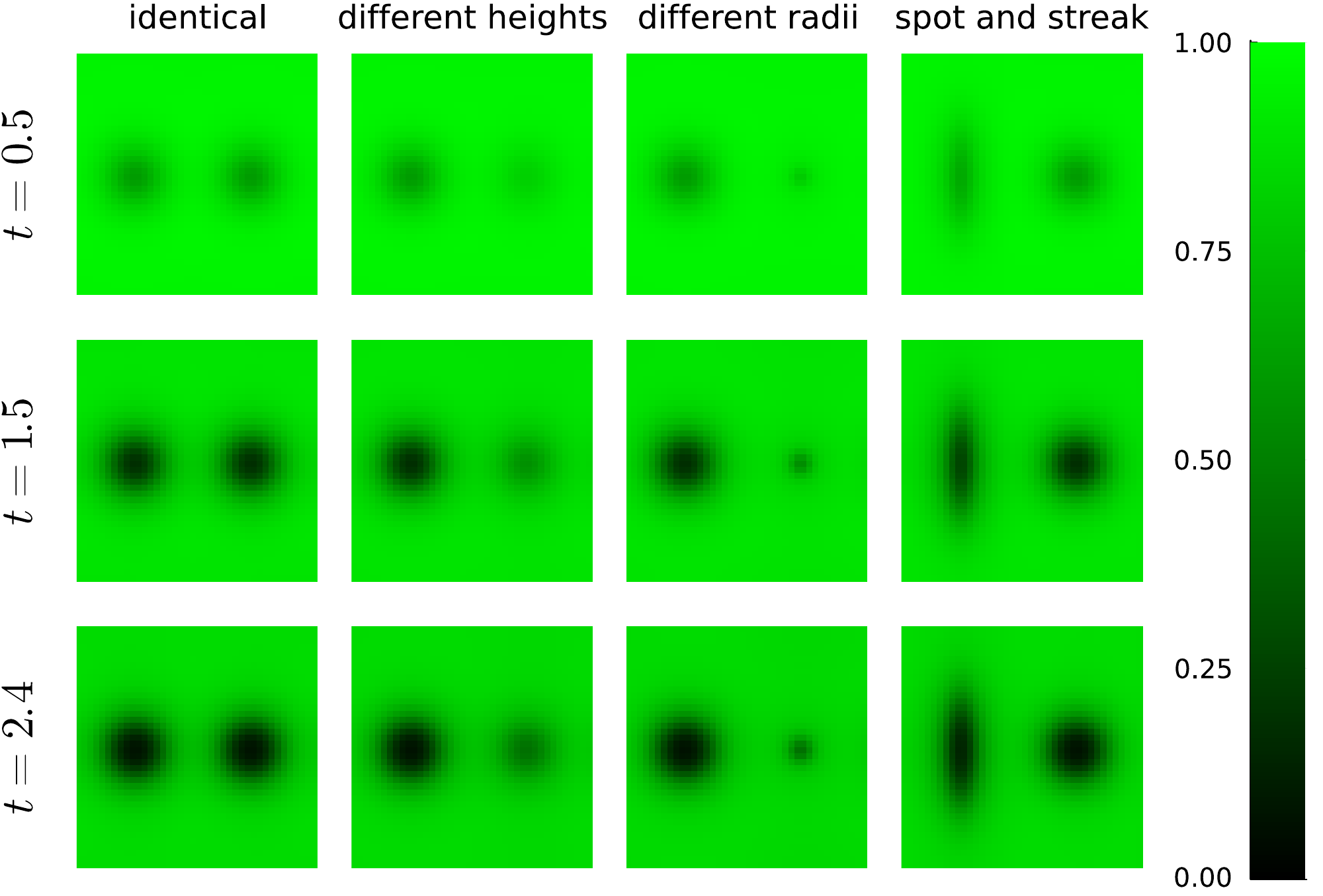}
    \caption{Snapshots of FL intensity for four cases of two well-separated evaporation regions; the parameters are given in \autoref{tab:two_spots_parameters}.}
    \label{fig:two_spots_separated_intensity}
\end{figure}

\autoref{fig:two_spots_separated_intensity} shows snapshots of the FL intensity for four different evaporation patterns with multiple thinning regions. The parameters of the evaporation function for these cases are given in \autoref{tab:two_spots_parameters}. From left to right, the columns represent: two identical spots; two same-sized spots with different peak evaporation rates; two spots with different radii but the same peak evaporation rates; and one spot and one thin ellipse. The TBUTs for all four simulations are, respectively, $2.6$, $2.58$, $2.49$, and $2.57$.

\begin{table}[tbp]
    \centering
    \begin{tabular}{c|ccc|ccc}
        Figure & \multicolumn{3}{c|}{Spot 1} & \multicolumn{3}{c}{Spot 2} \\ 
         & $a_1$ & $(x_1, y_1)$ & $(x_{w,1},y_{w,1})$ & $a_2$ & $(x_2, y_2)$ & $(x_{w,2},y_{w,2})$   \\
         \hline 
         \autoref{fig:two_spots_separated_intensity} (left) & $1$ & $(-1.5,0)$ & $(0.5, 0.5)$  & $1$ & $(1.5,0)$ & $(0.5, 0.5)$ \\
         \autoref{fig:two_spots_separated_intensity} (second to left) & $1$ & $(-1.5,0)$ & $(0.5, 0.5)$  & $0.5$ & $(1.5,0)$ & $(0.5, 0.5)$ \\
         \autoref{fig:two_spots_separated_intensity} (second to right) & $1$ & $(-1.5,0)$ & $(0.5, 0.5)$  & $1$ & $(1.5,0)$ & $(0.2, 0.2)$ \\
         \autoref{fig:two_spots_separated_intensity} (right) & $1$ & $(-1.5,0)$ & $(0.25, 1)$  & $1$ & $(1.5,0)$ & $(0.5, 0.5)$ \\  
    \end{tabular}
    \caption{Parameters in the evaporation function~\eqref{eq:Jm} for the simulations shown in \autoref{fig:two_spots_separated_intensity}. All have $K=2$ spots and baseline evaporation rate $v_b=0.1$.}
    \label{tab:two_spots_parameters}
\end{table}

 { The rest of this section focuses on three cases having evaporation patterns with identical circular spots that are separated by differing amounts.}  For the evaporation function $J(x,y)$  { in} \eqref{eq:Jm}, we set $a_k=1$ for both spots, $v_b=0.1${, and $x_{w,k} = y_{w,k} =0.5$ in} all three cases. We adjust the center positions, located at $\pm x_k$, using the three values $x_k = 1.5$, $0.8$, and $0.6$, respectively, and fix $y_k=0$.  The TBUT for each case is $2.6$, $2.3$, and $1.7$, respectively{, decreasing (as we show below) due to increased evaporative effects}.

\begin{figure}
    \centering
    \includegraphics[width=0.95\textwidth]{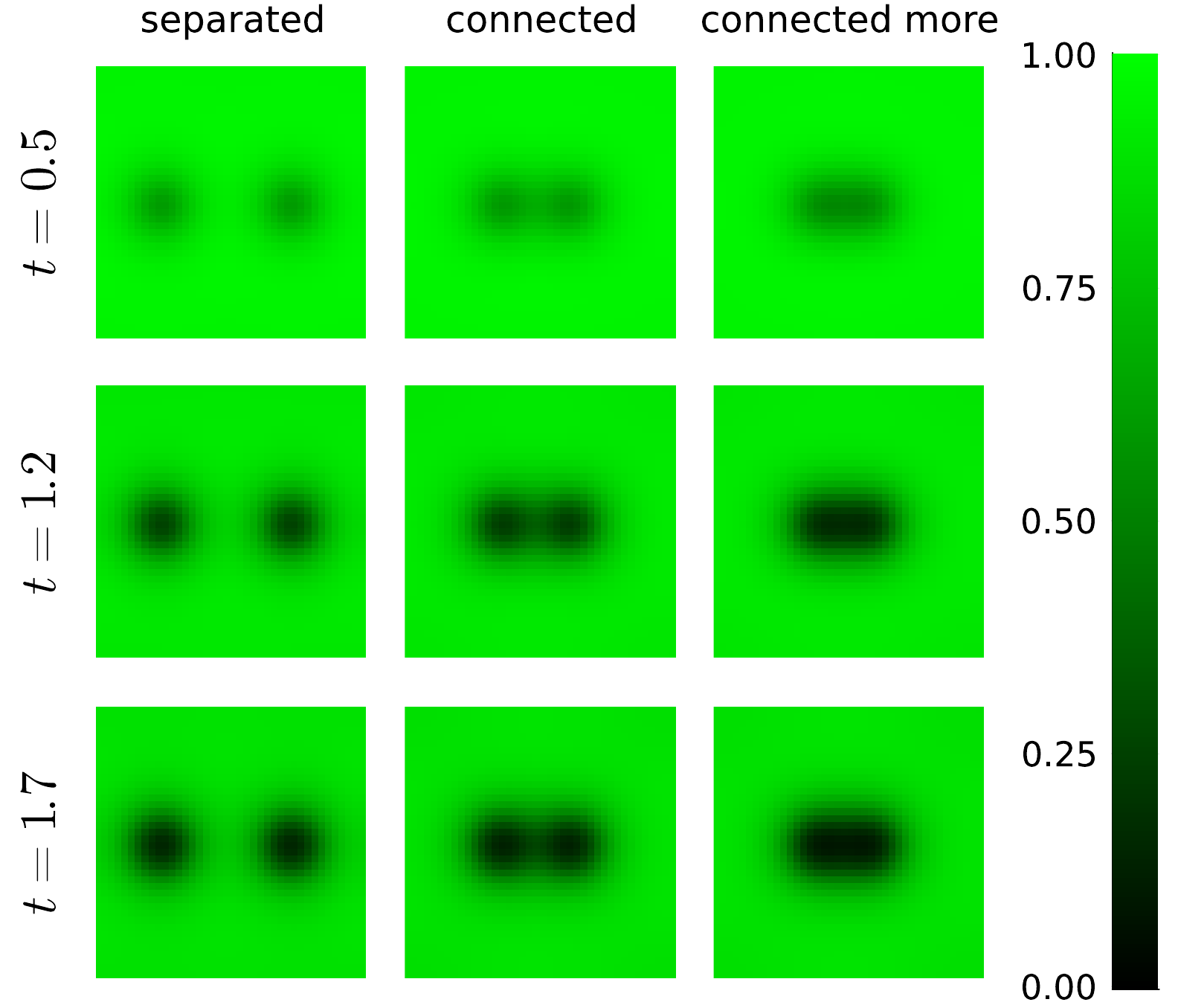}
    \caption{FL intensity contours at three different time levels for two spots from separated to connected. From left to right: $x_k=1.5$, $0.8$, and $0.6$.}
    \label{fig:two_spots_connected_intensity}
\end{figure}

{ \autoref{fig:two_spots_connected_intensity} shows the FL intensity at three time levels for each of these cases, with connectedness between the spots increasing from left to right.  Referring back to \autoref{fig:sptosk_fixed_xw_hcfi_center_and_connected_spots_center}, the right column} shows the central values {located at $\pm x_k$} of { the} dependent variables.   The central values for the largest  { separation} are very similar to the single-spot case, which is the blue line in the left column of { that figure,} though they occur slightly later because {of the larger TBUT.}  { As the separation decreases, the TBUT decreases as well.} The central values for the smaller $x_k$ are distinct from the central values from the elliptical cases in \autoref{fig:sptosk_fixed_xw_hcfi_center_and_connected_spots_center}. 
{\color{blue} 
Even though the TBUT is different, we observe that for $x_k=0.6$, which is the green dashed line in the second column, the osmolarity is very close to that of the streak case, indicating that the connected spots would behave similarly to a single streak or eccentric ellipse.}

\begin{figure}
    \centering
    \includegraphics[width=0.95\textwidth]{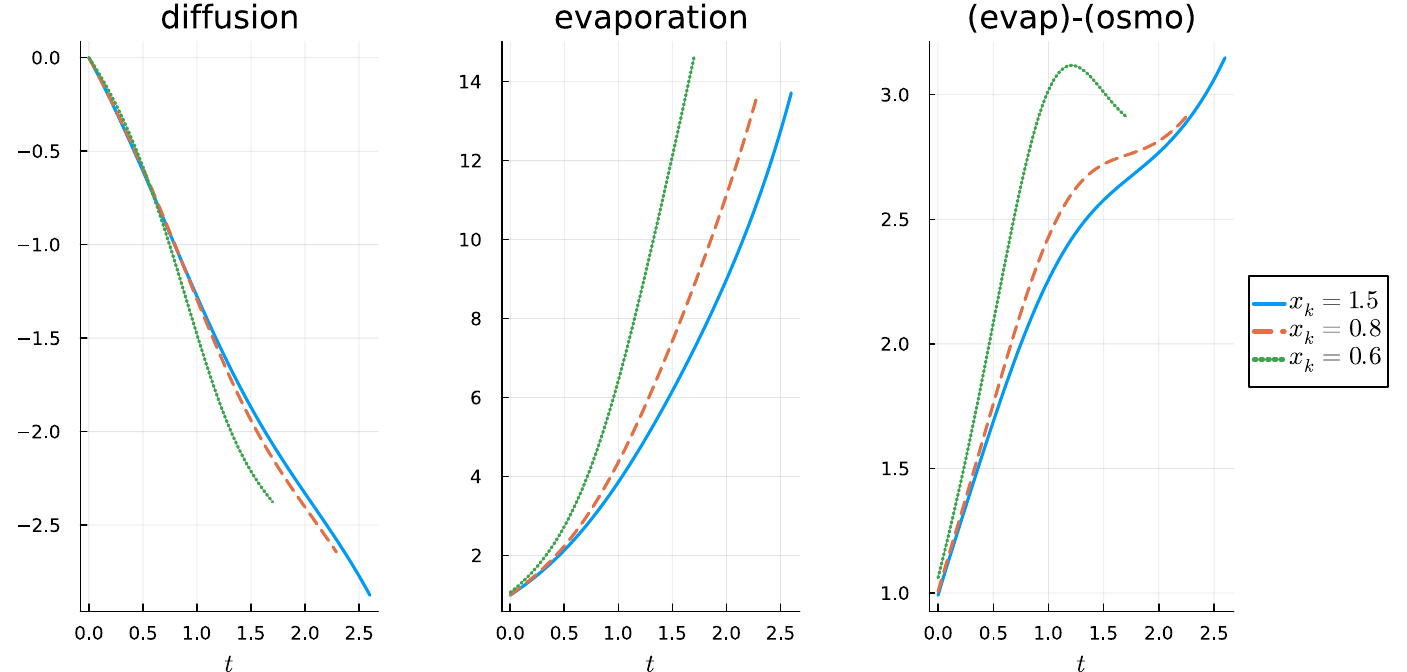}
        \caption{Dominant terms in the osmolarity equation \eqref{eq:termtypes}. The parameters for the evaporation function are the same as in the simulation for right column in \autoref{fig:sptosk_fixed_xw_hcfi_center_and_connected_spots_center}, \autoref{fig:two_connected_spots_hcfi_x_axis} and
        \autoref{fig:connected_osmo_xk}.}
    \label{fig:two_connected_spots_diffusion_center}
\end{figure}

\autoref{fig:two_connected_spots_diffusion_center} shows the dominant terms in the osmolarity equation \eqref{eq:termtypes}.  { We observe that as the separation distance decreases, the evaporation and diffusion increase, accounting for the faster TBU.} For $x_k=0.6$, the difference between evaporation and osmosis decreases in the end, showing that osmosis is increasing. 

\begin{figure}
    \centering
    \includegraphics[width=0.95\textwidth]{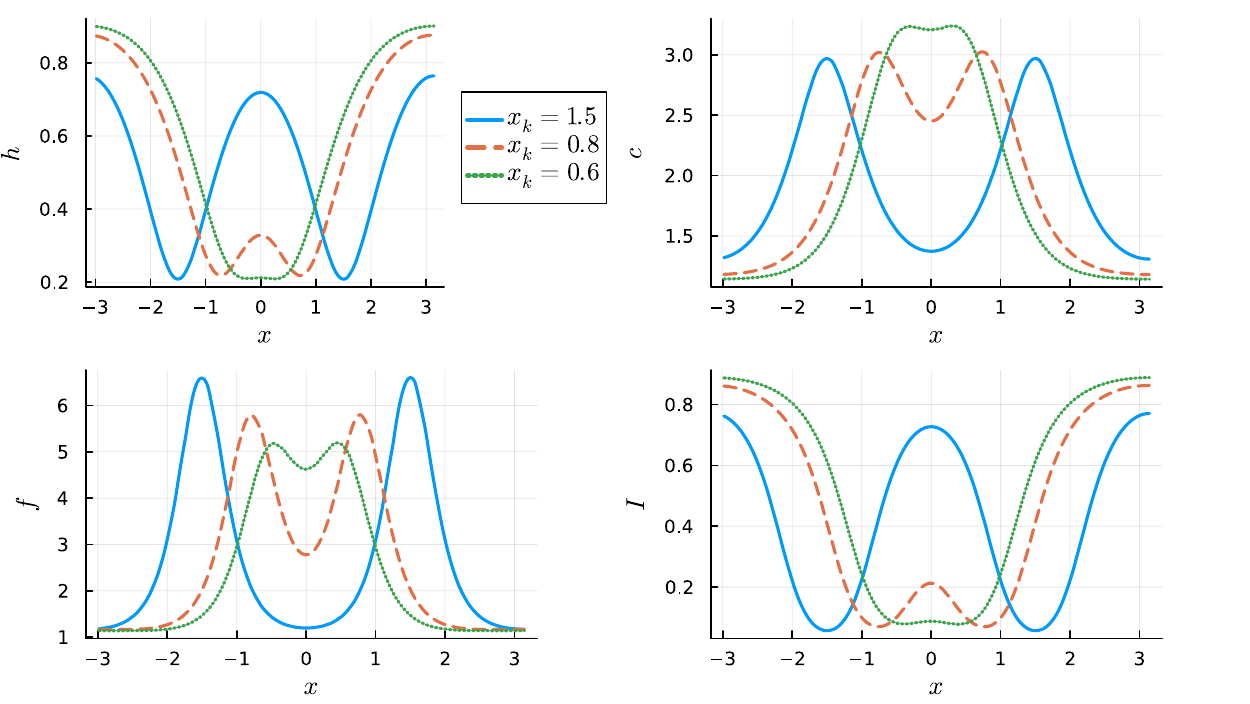}
        \caption{ { Values of $h,c,f,I$ along the separation axis at TBUT for identical spots with different separation half-distances $x_k$. }  }
    \label{fig:two_connected_spots_hcfi_x_axis}
\end{figure}

\autoref{fig:two_connected_spots_hcfi_x_axis} { shows the dependent variables along the $x$ axis (i.e., the separation axis) at TBUT. These plots are symmetric in $x$ because the spots are identical. The curves in the least-separated case $x_k=0.6$ suggest that the two spots have nearly fully merged.}

\begin{figure}
    \centering
    \includegraphics[width=0.8\textwidth]{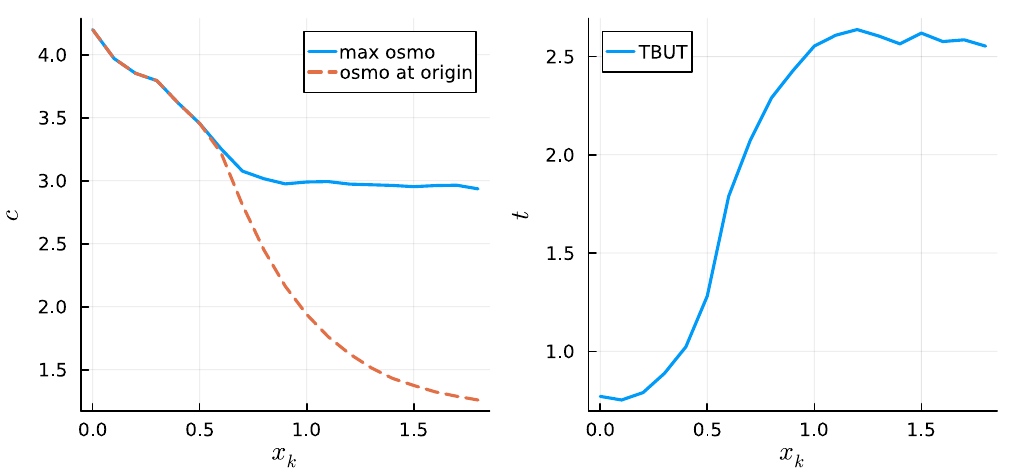}
        \caption{{ Left figure: Blue solid line shows the maximum osmolarity for two connected spots with varying $x_k$; Red dash line shows the osmolarity at the center of two connected spots. Right figure shows the TBUT with respect to $x_k$.}}
    \label{fig:connected_osmo_xk}
\end{figure}

{
Finally, \autoref{fig:connected_osmo_xk} shows in greater detail how the osmolarity at TBUT, and TBUT itself, vary with the separation distance. Up to a separation of $x_k=0.5$, which is the characteristic spot width $x_{w,k}$, the peak value of osmolarity at the spot centers practically coincides with the value at the origin, while the TBUT roughly doubles from the overlapped case. As the separation distance continues to increase, the peak osmolarity and TBUT level off at about $x_k=1$, suggesting that the spots no longer meaningfully interact beyond this distance. The asymptotic value of peak osmolarity is about 72\% of its full-overlap value, while the asymptotic TBUT is about 3 times the overlap value.
}             

{ \subsubsection{Performance of POD for connected spots}}

\begin{figure}
    \centering
    \includegraphics[width=0.95\textwidth]{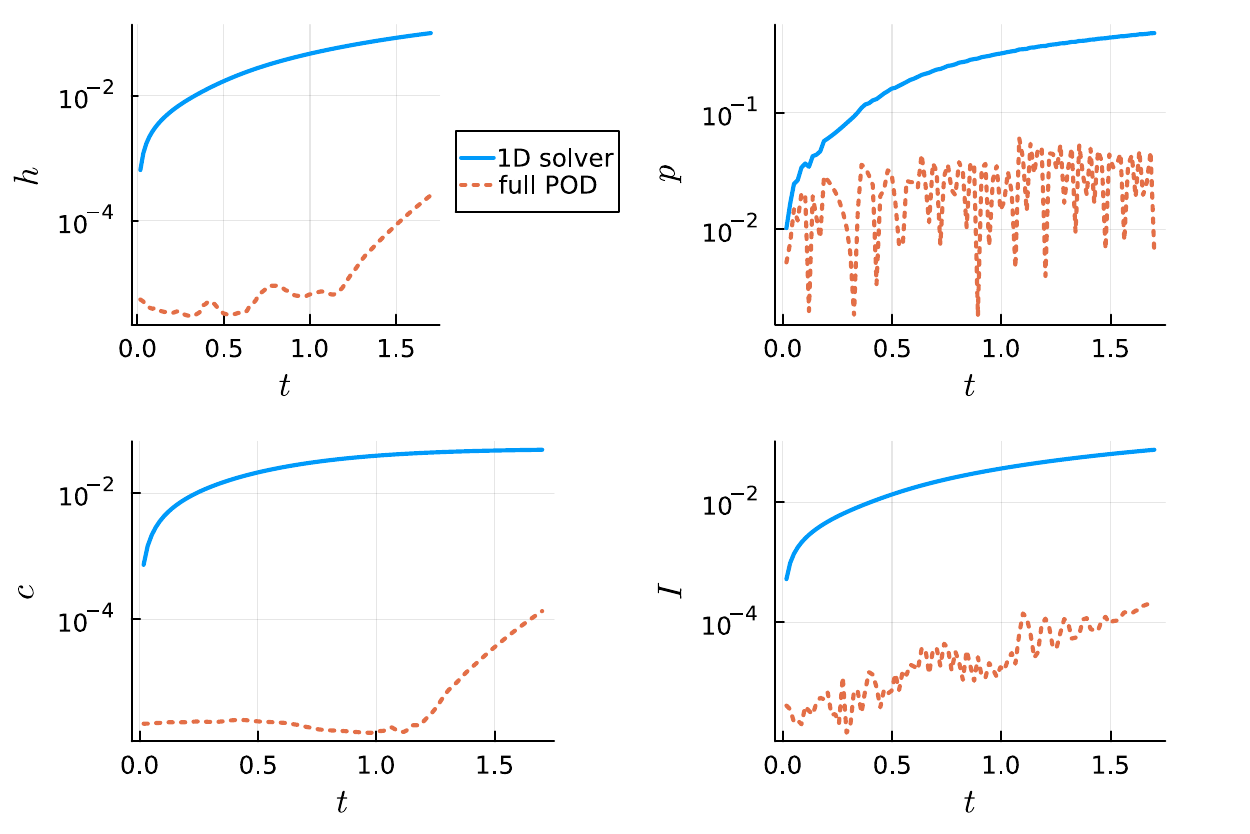}
    \caption{{Relative error for two variants of the POD dimension reduction method for two interacting evaporative spots.}}
    \label{fig:pod_err_connected_spots}
\end{figure}

{ 
For a single circular spot, we were able to use a 1D radial solver in place of the full 2D solution in order to quickly get numerical basis vectors for the POD method. With two well-separated spots, this technique can still be applied to the spots independently to construct the POD basis, but the error suffers as the separation decreases and the POD basis misses crucial information about the interaction between spots. \autoref{fig:pod_err_connected_spots} shows the relative error of the two POD variants  for the two connected spots in the third column in \autoref{fig:two_spots_connected_intensity}. The blue solid curve represents the POD solution obtained via 1D radial solver, and the red dashed curve represents the POD solution obtained via the 2D solver. The radial POD solution reaches essentially 100\% error well before TBUT, while the 2D POD stays below 10\% error. Although this observation is made purely from the standpoint of the performance of a numerical method, it provides further support for the conclusion that when two spots are close enough together, they behave very differently from two superimposed spots. 
}

\clearpage
\section{Conclusion}
\label{sec:conclusion}

{ In broad terms, there is a good resemblance between the time evolution of fluorescent intensity in our computational results, like those in \autoref{fig:sptosk_fixed_product_intensity} and \autoref{fig:sptosk_fixed_xw_intensity}, and in small patches within the experimental video data like those exhibited in \autoref{fig:flexample}. By varying the evaporation rate function, we can reproduce features seen in experiments, such as isolated spots and streaks as well as mixed shapes. We found, for instance, that once the aspect ratio of the ellipses for the evaporation function { reaches 2 (i.e., $y_w/x_w=2$)}, the TBU dynamics closely approximates a streak. In other words, when we have the relatively small distortion $y_w/x_w=2$, the central values of the solution from this elliptical distribution of evaporation are already a very good approximation to the solution for a streak ($y_w/x_w \to \infty$).} 

While our model does not account for every known effect within the tear film and ignores boundary effects near the eyelids, experience has shown~\cite{lukeParameterEstimation2020,lukeFittingSimplifiedModels2021,lukeParameterEstimationMixedMechanism2021,driscollFittingODEModels2023} that even limited models can be fit accurately to experimental FL data; in fact, our model can be seen as a direct 2D generalization of the 1D models in~\cite{lukeParameterEstimation2020}. Given the considerable modeling and computational challenges of complete 2D global behavior of the tear film, first investigating the utility of a restricted model is an attractive approach. 

Our next step will be to solve inverse problems for this model to get quantitative agreement with individual experiments, thus producing insights about thinning mechanisms and critical unobservable quantities such as local osmolarity. Since inverse problems necessitate multiple forward solutions and we intend to attempt fitting at many sites within individual trials over potentially hundreds of trials, computational efficiency of the forward solution is a major issue. This motivated our choice of Fourier discretization in space, which is the most efficient way to model smooth solutions in an interior domain. 

The inherent stiffness of the resulting DAE problem makes time-stepping difficult, and we have applied the POD method to project the system into far fewer dimensions to accelerate the solution process after a fixed time interval. The computational time for finding the reduced basis and solving the reduced model is insignificant compared to the full model. We have shown speedups of fourfold or more while keeping errors at a tolerable level.  Although we have not seen the POD method applied to systems with mass matrices elsewhere, that extension was straightforward, and there does not seem to be anything particular to this system that makes it especially challenging for POD. 

\section*{Declarations}
\textbf{Funding}. This work was supported by National Science Foundation grant 1909846. Any opinions, findings, and conclusions or recommendations expressed in this material are those of the authors and do not necessarily reflect the views of the National Science Foundation.

\vspace{1em} \noindent \textbf{Conflict of Interest}. The authors declare no competing interests.

\begin{appendices}
\section{Axisymmetric model}\label{secA1}

For the circular case, we use the axisymmetric coordinates $(r',z')$ to denote the position and $u'=(u',w')$ to denote the fluid velocity. The tear film is modeled as an incompressible Newtonian fluid on $0 < r' < R_0$ and $0<z'<h'(r',t')$. The scalings and non dimensional parameters are similar in \autoref{tab:physical} and \autoref{tab:parameters}. The only difference is that $r'=\ell r$. Braun et al. derived the system of equations on the domain $0<r<R_0$ \cite{braunTearFilm2018,braunDynamicsFunctionTear2015}.

\begin{align}
\partial_t h  + \partial_r (rh\overline{u}) &= -J + P_c(c-1), \label{eq:A1}\\
\overline{u} &=-\frac{h^2}{12}\partial_r p, \label{eq:A2}\\
p &=-\frac{1}{r}\partial_r(r\partial_r h), \label{eq:A3}\\
h(\partial_t c +\overline{u}\partial_r c ) &= \text{Pe}_c^{-1}\frac{1}{r}\partial_r (rh\partial_r c)+Jc-P_c(c-1)c, \label{eq:A4}\\
h(\partial_t f +\overline{u}\partial_r f ) &= \text{Pe}_f^{-1}\frac{1}{r}\partial_r (rh\partial_r f)+Jf-P_c(c-1)f, \label{eq:A5}
\end{align}
The evaporation function $J$ is

\begin{equation} \label{eq:AJ}
J(r) = v_b + (a-v_b) \exp\left[-(r/r_w)^2 / 2 \right]. 
\end{equation}
where $v_b$ is the ratio of $v_{\text{max}}$ over $v_{\text{min}}$, $r_w$ is the radius, and $a > v_b$ is the height of the peak.

\section{Full streak model}
The linear case model is solved on the Cartesian coordinates $-\pi < x < \pi$ and $ 0 < z < h(x,t)$. More details about derivation can be found in \cite{braunTearFilm2018,braunDynamicsFunctionTear2015}. Periodic boundary conditions are applied.  
\begin{align}
\partial_t h  + \partial_x (h\overline{u}) &= -J + P_c(c-1), \label{eq:B1}\\
\overline{u} &=-\frac{h^2}{12}\partial_x p, \label{eq:B2}\\
p &=-\partial_x^2 h, \label{eq:B3}\\
h(\partial_t c +\overline{u}\partial_x c ) &= \text{Pe}_c^{-1}\partial_x (h\partial_x c)+Jc-P_c(c-1)c, \label{eq:B4}\\
h(\partial_t f +\overline{u}\partial_x f ) &= \text{Pe}_f^{-1}\partial_x (h\partial_x f)+Jf-P_c(c-1)f, \label{eq:B5}
\end{align}
The evaporation function $J$ is

\begin{equation} \label{eq:AJ2}
J(x) = v_b + (a-v_b) \exp\left[-(x/x_w)^2 / 2 \right]. 
\end{equation}
The parameters are identical as in the spot case, and $r,r_w$ simply been replaced by $x,x_w$ here.

\end{appendices}

\bibliography{chenq_cite}

\end{document}